\documentclass[]{elsarticle}

\usepackage[title]{appendix}
\usepackage{amsmath}
\usepackage{mathtools}
\usepackage{amssymb}
\usepackage{upgreek}
\usepackage{fullpage}
\usepackage{color}
\usepackage{caption}
\usepackage{subcaption}
\usepackage{makecell}
\usepackage{mathrsfs}
\usepackage{breqn}
\usepackage{float}
\usepackage{tikz}
\usetikzlibrary{arrows,shapes,backgrounds,decorations.markings}

\usepackage[linesnumbered,ruled,vlined]{algorithm2e}
\SetKwInput{KwInput}{Input}                
\SetKwInput{KwOutput}{Output}              

\usepackage{float}
\usepackage{booktabs}
\usepackage{tabu} 
\usepackage{dcolumn}
\newcolumntype{A}{D{.}{.}{2.3}}
\usepackage{tabularx}
\usepackage[hidelinks]{hyperref}

\usepackage{amsthm}
\newcommand{\arch}{\mathrm{arch}}
\newcommand{\bb}{\mathrm{b}}
\newcommand{\comp}{\mathrm{comp}}
\newcommand{\fin}{\mathrm{fin}}
\newcommand{\port}{\mathrm{port}}
\newcommand{\train}{\mathrm{train}}
\newcommand{\test}{\mathrm{test}}
\newcommand{\rb}{\mathrm{rb}}
\newcommand{\sample}{\mathrm{sample}}
\newcommand{\sub}{\mathrm{sub}}
\newcommand{\mn}{\mathrm{min}}
\newcommand{\mx}{\mathrm{max}}
\newcommand{\lr}[1]{\left\| #1 \right\|}
\newcommand{\pr}[1]{\left( #1 \right)}

\DeclareMathOperator*{\argmin}{arg\,min}
\DeclareMathOperator*{\quadd}{\quad \quad}

\newcommand{\lib}{\what{\hspace{0pt}\mathcal{C}}}
\newcommand{\sys}{\mathcal{C}}
\newcommand{\pset}{\mathcal{P}}
\newcommand{\wpset}{\widehat{\mathcal{P}}}
\newcommand{\what}[1]{\widehat{#1}}
\newcommand{\wtilde}[1]{\widetilde{#1}}
\newcommand{\wg}{\widehat{\gamma}}
\newcommand{\wc}{\widehat{c}}
\renewcommand{\wp}{\widehat{p}}
\newcommand{\wrho}{\widehat{\rho}}
\renewcommand{\th}[1]{\tilde{\hat{#1}}}

\newcommand{\calO}{\mathcal{O}}

\newtheorem{theorem}{Theorem} 
\newtheorem{proposition}[theorem]{Proposition}
\newtheorem{corollary}[theorem]{Corollary}
\newtheorem{lemma}[theorem]{Lemma}
\newtheorem{remark}[theorem]{Remark}

\date{}

\begin{document}

\begin{frontmatter}

\title{A hyperreduced reduced basis element method for reduced-order modeling of component-based nonlinear systems\textsuperscript{*}}

\author[affl1,affl2]{Mehran Ebrahimi\corref{cor2}\fnref{fn1}}
\ead{m.ebrahimi@mail.utoronto.ca}
\author[affl1]{Masayuki Yano\fnref{fn2}}
\ead{masa.yano@utoronto.ca}

\affiliation[affl1]{
  organization={Institute for Aerospace Studies, University of Toronto},
  addressline={4925 Duffein Street},
  city={Toronto},
  postcode={M3H 5T6},
  state={Ontario},
  country={Canada}}
\affiliation[affl2]{
  organization={Autodesk Research},
  addressline={661 University Avenue},
  city={Toronto},
  postcode={M5G 1M1},
  state={Ontario},
  country={Canada}}
\cortext[cor1]{This is a preprint of an article published in \emph{Computer Methods in Applied Mechanics and Engineering}. The final authenticated version is available online at: \href{https://doi.org/10.1016/j.cma.2024.117254}{https://doi.org/10.1016/j.cma.2024.117254}.}
\cortext[cor2]{Corresponding author}
\fntext[fn1]{Graduate Student, Institute for Aerospace Studies, University of Toronto; Principal Research Scientist, Autodesk Research}
\fntext[fn2]{Associate Professor, Institute for Aerospace Studies, University of Toronto}

\begin{abstract}
\noindent
We introduce a hyperreduced reduced basis element method for model reduction of parametrized, component-based systems in continuum mechanics governed by nonlinear partial differential equations. In the offline phase, the method constructs, through a component-wise empirical training, a library of archetype components defined by a component-wise reduced basis and hyperreduced quadrature rules with varying hyperreduction fidelities. In the online phase, the method applies an online adaptive scheme informed by the Brezzi--Rappaz--Raviart theorem to select an appropriate hyperreduction fidelity for each component to meet the user-prescribed error tolerance at the system level. The method accommodates the rapid construction of hyperreduced models for large-scale component-based nonlinear systems and enables model reduction of problems with many continuous and topology-varying parameters. The efficacy of the method is demonstrated on a two-dimensional nonlinear thermal fin system that comprises up to 225 components and 68 independent parameters.
\end{abstract}


\end{frontmatter}

\section{Introduction}
\label{sec:introduction}
Many-query problems, which necessitate repeatedly solving parameterized partial differential equations (PDEs), arise commonly in various fields of computational science such as design optimization, uncertainty quantification, and control. For problems where the solution manifold is well approximated in a low-dimensional linear space, the reduced basis (RB) methods provide an effective approach to rapidly and reliably approximate the PDE solution for many different parameter values~\cite{rozza2008reduced, quarteroni2015reduced, hesthaven2016certified, benner2015survey}. RB methods achieve efficiency by separating the computation into offline (training) and online (deployment) phases. The former typically involves solutions of the high-fidelity problem (e.g., using finite element (FE) methods) for many training parameter values to generate \emph{solution snapshots}, the construction of an RB for the solution space, and, for nonlinear PDEs, \emph{hyperreduction}~\cite{patera2007reduced, rozza2008reduced, hesthaven2016certified}. Consequently, the offline phase is computationally demanding. Nonetheless, this initial high computational cost is warranted by the significant computational savings realized in the subsequent online phase, where the reduced problem is solved numerous times in the intended many-query application.

Despite their effectiveness, the applicability of standard (i.e., \emph{monodomain} or \emph{single-domain}) parametric RB methods is limited to the particular problem with continuous parametric variations for which the training is performed. For instance, even a slight topological change in the domain can render the trained reduced model inapplicable. In principle, a separate reduced model could be trained for each topological configuration; however, in practice, such retraining, at best, diminishes the utility of the reduced model and, at worst, is computationally intractable, especially for large-scale engineering systems that can take on many different topological configurations. Even when only parametric (and no topological) variations are considered, the standard RB methods can be restricted to problems with a small number of parameters due to the high training cost of exploring a high-dimensional space. 

To mitigate the aforementioned challenges, a variant of RB methods, called \emph{component-based} or \emph{multidomain} RB methods, have been developed~\cite{hurty1965dynamic, bourquin1992component, maday2002reduced}. The methods exploit the fact that many engineering structures---such as heat-exchangers, lattice structures, mechanical multi-component assemblies---consist of a large number of identical or similar components. The key ingredient of component-based RB methods is component-based training during the offline phase, whereby a library of interoperable \emph{archetype} components and their associated local RB is developed. Then, given a particular topological configuration in the online phase, copies of the archetype components in the library are \emph{instantiated}, and a global RB model for the whole system is formed by connecting the preconstructed local reduced models through their respective \emph{ports}.

Hitherto, several different variants of component-based RB methods have been developed. The reduced basis element (RBE) method~\cite{maday2002reduced,maday2004reduced,lovgren2006reduced} combines the ideas of domain-decomposition and RB methods. The method uses Lagrange multipliers to couple local, subdomain-wise reduced models in the online phase to form a global reduced model. The static condensation RBE (SCRBE) method~\cite{huynh2013static, huynh2013staticc} builds on the component mode synthesis~\cite{hurty1965dynamic, bourquin1992component} and the RB methods. The method decomposes the degrees of freedom (DoF) in each component into port and \emph{bubble} (interior) DoF. It then uses static condensation~\cite{wilson1974static} to form a Schur complement system with only port DoFs, and applies RB approximation in each component to reduce the computational cost of static condensation and to account for parametric variations. The port-reduced SCRBE method~\cite{eftang2013port, eftang2014port, smetana2015new} uses port-reduction techniques to further reduce the size of the Schur complement system and hence the computational cost. This is achieved by approximating the solution on global ports through the application of RB methods to port modes. SCRBE methods bear a close resemblance to multiscale RB methods~\cite{nguyen2008multiscale,boyaval2008reduced,kaulmann2011new,vidal2019multiscale,mcbane2021component,diercks2023multiscale}, which are applicable to structures composed of smaller-scale components with less heterogeneity relative to those in structures targeted for the SCRBE method. 

Component-based RB methods have been initially developed for linear or polynomial nonlinear problems with affine parameter dependence, which facilitate offline--online computational decomposition (without hyperreduction). Recently, these methods have been extended to general nonlinear and nonaffine problems. Methods for nonlinear problems can be broadly categorized into two groups based on the locality of nonlinearity. The first class of methods are designed for problems where the nonlinearity can be localized to small regions. Beiges et al.~\cite{baiges2013domain} decompose the physical domain and use a hybrid full-order/reduced-order model approach in the online phase to handle parameter configurations absent in the offline phase. Similarly, Ballani et al.~\cite{ballani2018component} decompose the physical domain into linear and nonlinear regions and apply the SCRBE method in the former part and high-fidelity model in the latter. Zhang et al.~\cite{zhang2019model} apply the same decomposition idea, but use Gaussian processes regression in nonlinear regions to construct a surrogate model. By construction, this class of methods is specialized for localizable nonlinearities and cannot treat globally nonlinear systems.

The second class of methods is designed for problems that exhibit nonlinearity everywhere in the domain. Hoang et al.~\cite{hoang2021domain} develop a domain-decomposition least-squares Petrov-Galerkin (DD-LSPG) method. This method constructs a separate reduced space for each subdomain and enforces interface continuity between the subdomains using a set of compatibility constraints in the LSPG method. Iollo et al.~\cite{iollo2023one} develop a component-based model reduction formulation for parametrized nonlinear elliptic PDEs that uses overlapping subdomains and an optimization-based reformulation. Smetana and Taddei~\cite{smetana2023localized} develop a multidomain RB method that uses the partition-of-unity concept and apply it to a two-dimensional nonlinear diffusion problem. Diaz et al.~\cite{diaz2024fast} integrate nonlinear approximation spaces, created through autoencoders, with domain-decomposition to facilitate reduced-order modeling of problems with slowly decaying Kolmogorov $n$-width. These methods for globally nonlinear problems, however, do not yet match the versatility and reliability offered by component-based RB methods for linear problems. First, the majority of these works consider multidomain systems that result from a decomposition of global system into partitions, and not interchangeable physical components in the sense of those in component-based RB methods for linear problems. Second, they do not provide a mechanism for quantitatively controlling the hyperreduction error at the system level during the online phase.

In this work, we propose a model reduction method that (i) can treat global nonlinearities, (ii) incorporates online-interchangeable physical components to provide topological and parametric online flexibility, and (iii) provides quantitative control of hyperreduction error. The contributions of the present work are fivefold:
\begin{enumerate}
\item We develop a hyperreduced RBE (HRBE) method, which (i) uses a library of archetype components to provide online topological and parametric flexibility of component-based RB methods \emph{and} (ii) can handle general parametrized nonlinear PDEs that exhibit global nonlinearities.
\item We extend the empirical quadrature procedure (EQP)~\cite{Patera_2017_EQP, yano2019lp} to component-wise offline training to enable a systematic construction of a library of hyperreduced components, each of which meets the specified residual tolerance.
\item We appeal to the Brezzi--Rappaz--Raviart (BRR) theory~\cite{caloz1997numerical} to develop an actionable solution error estimate for component-based nonlinear systems, which relates component-wise residuals due to hyperreduction to system-level solution error.
\item We develop an adaptive procedure, informed by the BRR error estimate, to construct a hyperreduced system from a library of hyperreduced components, such that the hyperreduction error in the online-assembled system meets the user-prescribed error tolerance in a solution norm for any topological and parametric configuration.
\item We demonstrate the efficacy of the proposed HRBE method using a nonlinear thermal fin system that comprises up to 225 instantiated components and 68 independent parameters.
\end{enumerate}

The remainder of the paper is organized as follows. Section~\ref{sec:def} presents the general form of the model problem considered in this study. Section~\ref{sec:nlrbe} introduces the HRBE method, providing the bubble--port decomposition, RB approximation, and hyperreduced RB approximation. Section~\ref{sec:training} introduces the component-based training procedure designed for RB construction and hyperreduction of the archetype components in the library. Section~\ref{sec:comppro} describes the computational procedures of offline and online phases. Section~\ref{sec:res} presents numerical results that validate and demonstrate the efficacy of the HRBE method. Finally, we conclude with a summary of the work and potential considerations for future work. 

\section{Parameterized nonlinear PDE model problem}
\label{sec:def}
As a prelude to developing our HRBE method, in this section, we introduce the general form of the considered parameterized nonlinear PDEs. We present both the physical and reference domain formulations, the latter of which is crucial to treat parameterized geometries using the HRBE method.

\subsection{Exact problem formulation}
\label{subsec:exact}
\begin{figure}[!bt]
\centering
    \begin{subfigure}{0.4\textwidth}
        \includegraphics[width=\textwidth]{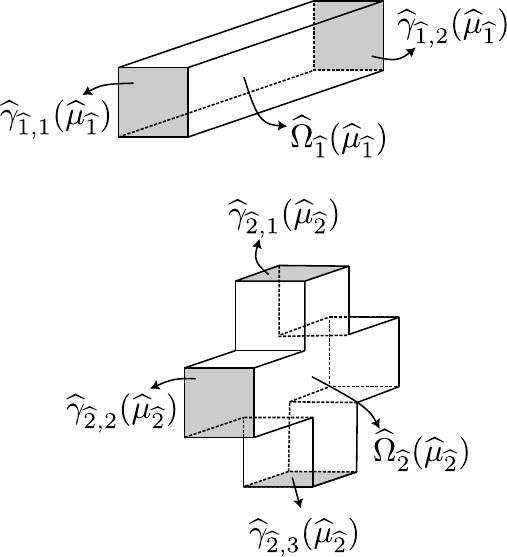}
        \caption{components}\label{fig:sysexamplea}
    \end{subfigure}
    \hspace{5em}
    \begin{subfigure}{0.45\textwidth}
        \includegraphics[width=0.9\textwidth]{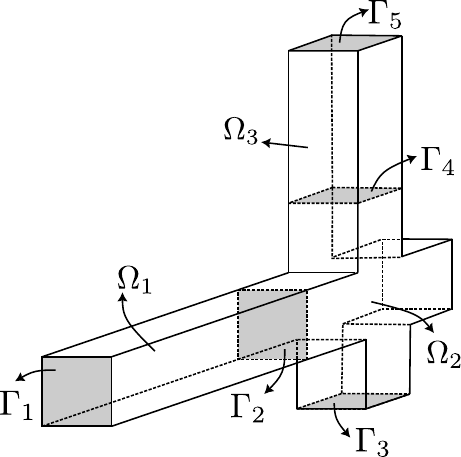}
        \caption{assembled system}\label{fig:sysexampleb}
    \end{subfigure}
    \caption{(a) Top: an archetype component with two local ports and $\wpset_1 = \{1,2 \}$, Bottom: an archetype component with three local ports and $\wpset_2 = \{1,2,3 \}$; (b) A system with $N_\comp = 3$ instantiated components and $N_\port = 5$ global ports. In this system, $M(1) = \what{1}$, $M(2) = \what{2}$, $M(3) = \what{1}$, and $\pset = \{1,\cdots,5 \}$. Also, $\Omega_1 = \mathcal{G}_1(\what{\Omega}_{M(1)};\mu_1)$, $\Omega_2 = \mathcal{G}_2(\what{\Omega}_{M(2)};\mu_2)$, $\Omega_3 = \mathcal{G}_3(\what{\Omega}_{M(3)};\mu_3)$, $\Gamma_1 = \mathcal{G}_1(\what{\gamma}_{M(1), {1}};\mu_1)$, $\Gamma_2 = \mathcal{G}_1(\what{\gamma}_{M(1), {2}};\mu_1) = \mathcal{G}_2(\what{\gamma}_{M(2), {2}};\mu_2)$, $\Gamma_3 = \mathcal{G}_2(\what{\gamma}_{M(2), {3}};\mu_2)$, $\Gamma_4 = \mathcal{G}_2(\what{\gamma}_{M(2), {1}};\mu_2) = \mathcal{G}_3(\what{\gamma}_{M(3), {1}};\mu_3)$, and $\Gamma_5 = \mathcal{G}_3(\what{\gamma}_{M(3), {2}};\mu_3)$.}
    \label{fig:sysexample}
\end{figure}

We first introduce geometric and topological entities associated with archetype components. We define $\lib$ as a \emph{library} of $N_\arch$ \emph{parameterized archetype components}. For each archetype component $\what{c} \in \lib$, we introduce $\what{\Omega}_{\what{c}} \subset \mathbb{R}^d$, ${\what{\mathcal{D}}}_{\what{c}} \subset \mathbb{R}^{{n}_{\what{c}}}$, and ${\what{\mu}}_{\what{c}} \in {\what{\mathcal{D}}}_{\what{c} }$ as, respectively, its bounded $d$-dimensional reference spatial domain, bounded ${n}_{\what{c}}$-dimensional parameter domain, and ${n}_{\what{c}}$-dimensional parameter tuple specifying its reference parameter values. Each archetype component $\what{c}$ has $n_{\what{c}}^{\gamma}$ disjoint \emph{local ports} whose domains are $\what{\gamma}_{\what{c},\what{p}} \subseteq \partial \what{\Omega}_{\what{c}}$, $\what{p} \in {\mathcal{P}}_\wc \equiv \{ 1, \cdots, n_{\wc}^{\gamma} \}$, where $\partial \what{\Omega}_{\what{c}}$ is the boundary of $\what{\Omega}_{\what{c}}$. We assume the boundary of all components are Lipschitz continuous and all ports of an archetype component are mutually separated by a boundary surface. Figure~\ref{fig:sysexamplea} shows these definitions for two archetype components.

We next introduce geometric and topological entities associated with an assembled system. We define $\sys$ as a set of $N_\comp$ \emph{instantiated components} composing a system. Each instantiated component is generated from an archetype component in the library through a (parameterized) geometric mapping. The components in the system are connected together through their local ports, thereby creating $N_\port$ \emph{global ports}. The geometric mappings must guarantee compatibility of the ports. We assume a global port can be shared by at most two instantiated components. A local port residing on the system boundary also forms a global port, where the essential boundary conditions at the system level are imposed. We introduce $\Omega_c \subset \mathbb{R}^d$ as the physical domain of the instantiated component $c \in \sys$, and $\Gamma_p$, $p \in \pset \equiv \{ 1, \cdots, {N_\port} \}$, as the physical domain of the $p$-th global port in the system. We introduce, for each instantiated component $c$, the parameter tuple $\mu_c \in \mathcal{D}_c \equiv \widehat{\mathcal{D}}_{{M}(c)}$, where $M: \sys \rightarrow \lib$ is a map from the instantiated components to their corresponding archetype components. The parameterized geometric mappings relating the archetype and instantiated component domains are $\mathcal{G}_c: \widehat{\Omega}_{{M}(c)} \times \mathcal{D}_{c} \rightarrow \Omega_c$ such that $\Omega_c = \mathcal{G}_c (\widehat{\Omega}_{{M}(c)}; \mu_c)$. The physical domain of the $p$-th local port of $c \in \sys$ is given by $\gamma_{c,p} \equiv \mathcal{G}_c (\widehat{\gamma}_{{M}(c),p}; \mu_c) \: \: \forall p \in {\pset}_{M(c)}$. The mapping $\mathcal{G}_c(\cdot;\mu_c)$ depends only on the geometric parameters in $\mu_c$. Figure~\ref{fig:sysexampleb} shows an example of a three-component system.

We now define function spaces associated with archetype and instantiated components. For the archetype component ${\what{c}} \in \lib$, we introduce a Hilbert space $\widehat{\mathcal{V}}_{{\what{c}}} \subset H^1(\widehat{\Omega}_{\what{c}})$ endowed with an inner product $(\cdot, \cdot)_{\widehat{\mathcal{V}}_{{\what{c}}}}$ and the associated induced norm $|| \cdot ||_{\widehat{\mathcal{V}}_{{\what{c}}}} \equiv \sqrt{(\cdot, \cdot)_{\widehat{\mathcal{V}}_{{\what{c}}}}}$, which is equivalent to the $H^1(\widehat{\Omega}_{\what{c}})$-norm. For $c \in \sys$, we introduce the geometric-parameter-dependent mapped space $\mathcal{V}_{c} \equiv \left\{v = \widehat{v} \circ \mathcal{G}_c^{-1}(\cdot; \mu_c) \Big| \: \widehat{v} \in \widehat{\mathcal{V}}_{{M}(c)} \right\}$ and the associated inner product $(\cdot, \cdot)_{\mathcal{V}_c}$ and induced norm $\|\cdot \|_{\mathcal{V}_c} \equiv \sqrt{(\cdot, \cdot)_{\mathcal{V}_c}}$.

We now present a domain-decomposed formulation of the system-level model problem in terms of its components. We introduce the system's physical domain $\Omega$ such that $\overline \Omega = \cup_{c \in \sys} \overline \Omega_c$. We also define $\Gamma_D$ and $\Gamma_N$, respectively, as the nonempty Dirichlet and Neumann boundaries of $\Omega$ such that $\partial \Omega = \overline{\Gamma}_D \cup \overline{\Gamma}_N$ and $\Gamma_D \cap \Gamma_N = \emptyset$. The Dirichlet boundary is composed of nonshared local ports of the instantiated components in the system. To simplify the presentation, we assume homogeneous boundary conditions everywhere. We further introduce the system parameter domain $\mathcal{D} \equiv \oplus_{c \in \sys} \mathcal{D}_{c}$ and parameter tuple $\mu \equiv (\mu_1, \cdots, \mu_{N_\comp}) \in \mathcal{D}$. Additionally, for each instantiated component $c \in \sys$, we introduce the physical-domain residual ${R}_c : {\mathcal{V}}_c \times {\mathcal{V}}_c \times {\mathcal{D}}_c \rightarrow \mathbb{R}$ as
\begin{equation}\nonumber
\begin{aligned}
{R}_c (w,v;\mu) &= \int_{{\Omega}_c} {r}_c(w,v;{x},\mu)  \: d{x} & &\quad \forall w,v \in {\mathcal{V}}_c, \: \forall \mu \in {\mathcal{D}}_c,
\end{aligned}
\end{equation}
where ${r}_c: {\mathcal{V}}_c \times {\mathcal{V}}_c \times {\Omega}_c \times {\mathcal{D}}_c \rightarrow \mathbb{R}$ is the physical-domain integrand, which is linear in its second argument but is in general nonlinear in its first argument. The \emph{exact} nonlinear model problem in its weak form is as follows: given $\mu = (\mu_c)_{c \in \sys} \in \mathcal{D}$, find $u(\mu) \in\ \mathcal{V}$ such that
\begin{equation}\label{eq:exact}
\begin{aligned}
R(u(\mu), v;\mu) \equiv \sum_{c \in \sys}  {R}_{c} \left( u(\mu) \big|_{\Omega_c}, v \big|_{\Omega_c}; \mu_c  \right) = 0 \quad \forall v \in \mathcal{V},
\end{aligned}
\end{equation}
where $\mathcal{V} = \left\{ v \in H^1 (\Omega)  \Big| v |_{\Gamma_D} = 0 \right\}$. Problems that involve boundary integrals due to nonhomogeneous boundary conditions can be readily treated with minor modifications. We assume the problem is well-posed for all $\mu \in \mathcal{D}$. Given the solution field $u(\mu) \in \mathcal{V}$, we evaluate a scalar output (i.e., quantity of interest) $F(u(\mu);\mu) \in \mathbb{R}$ at the system level, where
\begin{equation}\nonumber
\begin{aligned}
{F} (w;\mu) &= \sum_{c \in \sys} {F}_c (w \big|_{\Omega_c};\mu_c) \quad \forall w \in \mathcal{V}, \forall \mu \in \mathcal{D};
\end{aligned}
\end{equation}
here, ${F}_c : {\mathcal{V}}_c \times {\mathcal{D}}_c \rightarrow \mathbb{R}$ is the physical-domain output functional for the instantiated component $c \in \sys$ given by
\begin{equation}\nonumber
    \begin{aligned}
    {F}_c (w;\mu) &= \int_{{\Omega}_c} {f}_c(w;{x},\mu)  \: d{x} & &\quad \forall w \in {\mathcal{V}}_c, \: \forall \mu \in {\mathcal{D}}_c,
\end{aligned}
\end{equation}
where ${f}_c: {\mathcal{V}}_c \times {\Omega}_c \times {\mathcal{D}}_c \rightarrow \mathbb{R}$ is the physical-domain integrand.

To handle parameterized geometries using the HRBE method presented in the next section, we need to formulate the system-level residual and output forms in the reference-domain of the components. As such, for each archetype component $\wc \in \lib$, we introduce reference-domain residual $\widehat{R}_\wc : \widehat{\mathcal{V}}_\wc \times \widehat{\mathcal{V}}_\wc \times \widehat{\mathcal{D}}_\wc \rightarrow \mathbb{R}$ and output functional $\widehat{F}_\wc : \widehat{\mathcal{V}}_\wc \times \widehat{\mathcal{D}}_\wc \rightarrow \mathbb{R}$ given by
\begin{equation}\nonumber
\begin{aligned}
\widehat{R}_\wc (w,v;\mu) &= \int_{\widehat{\Omega}_\wc} \widehat{r}_\wc (w,v;\widehat{x},\mu)  \: d\widehat{x} & &\quad \forall w,v \in \widehat{\mathcal{V}}_\wc, \: \forall \mu \in \widehat{\mathcal{D}}_\wc, \\
\widehat{F}_\wc (w;\mu) &= \int_{\widehat{\Omega}_\wc} \widehat{f}_\wc (w;\widehat{x},\mu)  \: d\widehat{x} & &\quad \forall w \in \widehat{\mathcal{V}}_\wc, \: \forall \mu \in \widehat{\mathcal{D}}_\wc,
\end{aligned}
\end{equation}
where $\widehat{r}_\wc: \widehat{\mathcal{V}}_\wc \times \widehat{\mathcal{V}}_\wc \times \widehat{\Omega}_\wc \times \widehat{\mathcal{D}}_\wc \rightarrow \mathbb{R}$ and $\widehat{f}_\wc: \widehat{\mathcal{V}}_\wc \times \widehat{\Omega}_\wc \times \widehat{\mathcal{D}}_\wc \rightarrow \mathbb{R}$ are the reference-domain integrands that satisfy
\begin{equation}\nonumber
    \begin{aligned}
    {r}_c(w,v;{x},\mu_c) &= \what{r}_{M(c)}(w \circ \mathcal{G}_c(\cdot; \mu_c),v \circ \mathcal{G}_c(\cdot; \mu_c); \mathcal{G}_c^{-1}(x; \mu_c),\mu_c)  \: \Big| \mathcal{J}_c(\mathcal{G}_c^{-1}(x; \mu_c); \mu_c) \Big|^{-1}, \\
    {f}_c(w;{x},\mu_c) &= \what{f}_{M(c)}(w \circ \mathcal{G}_c(\cdot; \mu_c); \mathcal{G}_c^{-1}(x; \mu_c),\mu_c)  \: \Big| \mathcal{J}_c(\mathcal{G}_c^{-1}(x; \mu_c); \mu_c) \Big|^{-1} \\
\end{aligned}
\end{equation}
for all $c \in \sys$, $w,v \in \mathcal{V}_c$, $x \in \Omega_c$, and $\mu_c \in \mathcal{D}_c$; here, $\mathcal{J}_c(\cdot; \mu_c)$ is the Jacobian of $\mathcal{G}_c(\cdot; \mu_c)$, and $\Big| \mathcal{J}_c(\cdot; \mu_c) \Big|$ is its determinant. We can now express the system-level residual and output forms in terms of the archetype component reference-domain forms as
\begin{equation}\nonumber
\begin{aligned}
R(w, v;\mu) &= \sum_{c \in \sys}  \widehat{R}_{M(c)} \left( w \big|_{\Omega_c} \circ \mathcal{G}_c(\cdot; \mu_c), v \big|_{\Omega_c} \circ \mathcal{G}_c(\cdot; \mu_c); \mu_c  \right), \\
{F} (w;\mu) &= \sum_{c \in \sys} \what{F}_{M(c)} (w \big|_{\Omega_c} \circ \mathcal{G}_c(\cdot; \mu_c);\mu_c)
\end{aligned}
\end{equation}
for any $w,v \in \mathcal{V}$ and $\mu \in \mathcal{D}$.

\subsection{Truth problem formulation}
\label{subsec:truth}
As is often the case, the exact problem~\eqref{eq:exact} cannot be solved analytically. Instead, we appeal to the \emph{truth} problem associated with a FE method to approximate the solution. This solution is taken as the computable ground truth. We present the truth problem formulation in terms of bubble and port functions to facilitate the development of the HRBE method described in Section~\ref{sec:nlrbe}.

\subsubsection{Bubble--port decomposition of functions}
We first define approximation spaces associated with archetype components. For each archetype component $\wc \in \lib$, we introduce an $\mathcal{N}_{\what{c}}$-dimensional \emph{truth} FE space $\widehat{\mathcal{V}}_{h,{\what{c}}} \equiv \left\{ v \in \widehat{\mathcal{V}}_{{\what{c}}} \Big| \: v|_\kappa \in \mathbb{P}^n(\kappa) \: \: \forall \kappa \in \mathcal{T}_{h,{\what{c}}} \right\} \subset \widehat{\mathcal{V}}_{{\what{c}}}$, where $\mathcal{T}_{h,{\what{c}}}$ is a tessellation of $\widehat{\Omega}_{\what{c}}$ formed by a set of nonoverlapping, conforming elements $\{ \kappa \}$, and $\mathbb{P}^n(\kappa)$ is the space of degree-$n$ polynomials over each element $\kappa$. We also introduce $\widehat{\mathcal{V}}_{h,{\what{c}}}^\bb \equiv \left\{ v \in \widehat{\mathcal{V}}_{h,{\what{c}}} \Big| \: v \big|_{\wg_{\wc,\wp}} = 0, \ \forall \wp \in {\pset}_{\what{c}} \right\}$ as the $\mathcal{N}_{\what{c}}^\bb$-dimensional bubble FE space of the archetype component ${\what{c}}$. We additionally introduce $\mathcal{N}_{{\what{c}}}^{\wp}$-dimensional port FE space $\widehat{\mathcal{X}}_{h,{\what{c}}}^{\wp}$ of the $\wp$-th port of the archetype component ${\what{c}}$ as the restriction of $\widehat{\mathcal{V}}_{h,{\what{c}}}$ to the port domain $\wg_{\wc,\wp}$; i.e., $\widehat{\mathcal{X}}_{h,{\what{c}}}^{\wp} \equiv \widehat{\mathcal{V}}_{h,{\what{c}}}|_{\wg_{\wc,\wp}}$, $\wp \in \pset_\wc$. We note that $\mathcal{N}_{\what{c}} = \mathcal{N}_{\what{c}}^\bb + \sum_{\wp \in {\pset}_{\what{c}}} \mathcal{N}_{{\what{c}}}^\wp$. 

We now define basis functions for the ports of each archetype component $\wc \in \lib$. We introduce for the $\wp$-th local port of $\wc$ eigenpairs $( \widehat{\tau}_{\wc,i}^\wp \in \widehat{\mathcal{X}}_{h,\wc}^\wp, \lambda_{\wc,i}^\wp \in \mathbb{R} )_{i=1}^{\mathcal{N}_{\wc}^\wp}$
such that
\begin{equation}\nonumber
\begin{aligned}
\int_{\widehat{\gamma}_{\wc,\wp}} \nabla \widehat{\tau}_{\wc,i}^\wp \cdot \nabla \what{y} \: ds &= \lambda_{\wc,i}^\wp \int_{\widehat{\gamma}_{\wc,\wp}} \widehat{\tau}_{\wc,i}^\wp \: \what{y} \: ds \quad \forall \what{y} \in \widehat{\mathcal{X}}_{h,\wc}^\wp, \\
\left\| \widehat{\tau}_{\wc,i}^\wp \right\|_{L^2(\widehat{\gamma}_{\wc,\wp})} &= 1.
\end{aligned}
\end{equation}
We note that $\widehat{\mathcal{X}}_{h,{\what{c}}}^{\wp} = \text{span} \{ \widehat{\tau}_{\wc,i}^\wp \}_{i=1}^{\mathcal{N}_{\wc}^\wp}$. We then elliptically lift these basis functions to the interior of $\wc$ to find $\{ \widehat{\psi}_{\wc,i}^\wp \in \widehat{\mathcal{V}}_{h,\wc} \}_{i=1}^{\mathcal{N}_{\wc}^\wp}$ by solving
\begin{equation}\nonumber
\begin{alignedat}{3}
\int_{\widehat{\Omega}_\wc} \nabla \widehat{\psi}_{\wc,i}^\wp \cdot \nabla v \:  d\widehat{x} &= 0  \quad \quad && \forall v \in \widehat{\mathcal{V}}_{h,\wc}^\bb, \\
\widehat{\psi}_{\wc,i}^\wp &= \widehat{\tau}_{\wc,i}^\wp \quad \quad && \text{on } \widehat{\gamma}_{\wc,\wp}, \\
\widehat{\psi}_{\wc,i}^\wp &= 0 \quad \quad && \text{on } \widehat{\gamma}_{\wc,\wp\:'} \quad \forall \wp\:' \neq \wp.
\end{alignedat}
\end{equation}
We define $\what{\mathcal{V}}_{h,\wc}^\gamma$ as the $\mathcal{N}_\wc^\gamma$-dimensional FE space spanned by $\{ \{ \what{\psi}_{\wc,i}^{\wp} \}_{i=1}^{\mathcal{N}_{\wc}^\wp}\}_{\wp \in \pset_{\wc}}$, where $\mathcal{N}_\wc^\gamma = \sum_{\wp \in \pset_\wc} \mathcal{N}^{\wp}_\wc$.

We now present the bubble--port decomposition of functions defined on each instantiated component $c \in \sys$. We introduce the geometric-parameter-dependent mapped full-component, bubble, port(-trace), and port-lifted FE spaces:
\begin{equation}
    \begin{aligned}
        \mathcal{V}_{h,c} &\equiv \left\{v = \widehat{v} \circ \mathcal{G}_c^{-1}(\cdot; \mu_c) \Big| \: \widehat{v} \in \widehat{\mathcal{V}}_{h, {M}(c)} \right\} \subset \mathcal{V}_{c}, \nonumber \\
        \mathcal{V}_{h,c}^\bb &\equiv \Big\{v = \widehat{v} \circ \mathcal{G}_c^{-1}(\cdot; \mu_c) \Big|\: \widehat{v} \in \widehat{\mathcal{V}}_{h, {M}(c)}^\bb \Big\} \subset \mathcal{V}_{h,c}, \nonumber \\
        \mathcal{X}_{h,c}^p &\equiv \Big\{v = \widehat{v} \circ \mathcal{G}_c^{-1}(\cdot; \mu_c) \Big|\: \widehat{v} \in \widehat{\mathcal{X}}_{h, {M}(c)}^p \Big\} \quad \forall p \in \pset_{M(c)}, \nonumber \\
        \mathcal{V}_{h,c}^\gamma &\equiv \Big\{v = \widehat{v} \circ \mathcal{G}_c^{-1}(\cdot; \mu_c) \Big|\: \widehat{v} \in \widehat{\mathcal{V}}_{h, {M}(c)}^\gamma \Big\} = \text{span} \Big\{ \Big\{ \psi_{c,i}^p \equiv \widehat{\psi}_{M(c),i}^p \circ \mathcal{G}_c^{-1}(\cdot;\mu_c) \Big\}_{i=1}^{\mathcal{N}_{M(c)}^p} \Big\}_{p \in \pset_{M(c)}}.
    \end{aligned}
\end{equation}
Subsequently, any $v_{h,c} \in \mathcal{V}_{h,c}$ can be decomposed as
\begin{equation}\label{eq:compsol}
    v_{h,c} = v^\bb_{h,c} + v^\gamma_{h,c},
\end{equation}
where $v^\bb_{h,c} \in \mathcal{V}_{h,c}^\bb$ is the \emph{bubble part} of $v_{h,c}$ and $v^\gamma_{h,c} \in \mathcal{V}_{h,c}^\gamma$ is its \emph{port part} given by
\begin{equation}\nonumber
    v^\gamma_{h,c} = \sum_{p\in \pset_{M(c)}} v^p_{h,c} = \sum_{p\in \pset_{M(c)}} \sum_{i=1}^{\mathcal{N}^p_{M(c)}} \mathbf{v}^p_{h,c, i} \psi^p_{c,i},
\end{equation}
in which $\{ \mathbf{v}^p_{h,c, i} \}_{i=1}^{\mathcal{N}^p_{M(c)}}$ are the \emph{generalized coordinates} of $v^p_{h,c} \in \mathcal{V}_{h,c}^p \equiv \text{span} \Big\{ \psi_{c,i}^p \Big\}_{i=1}^{\mathcal{N}_{M(c)}^p}$ $\forall c \in \sys$ and $\forall p \in \pset_{M(c)}$. (Throughout this work, we denote the generalized coordinates of any function $y$ in an $N$-dimensional linear space $\mathcal{Y}$ with a basis $\{ \Phi_i \}_{i=1}^N$ by a boldface letter $\mathbf{y}$ such that $\mathbf{y} \equiv [\mathbf{y}_1,\cdots,\mathbf{y}_N]^T$ satisfies $y = \sum_{i=1}^N \mathbf{y}_i \Phi_i$.)

At the system level, we assume conformity of the local ports connected together. Thus, for the $p$-th global port, $p \in \pset$, shared by $l$-th port of $c \in \sys$ and $l'$-th port of $c' \in \sys$, we have $\Gamma_p = \gamma_{c,l} = \gamma_{c',l'}$, $\mathcal{X}_{h,p} \equiv \mathcal{X}_{h,c}^l = \mathcal{X}_{h,c'}^{l'}$, $\mathcal{N}^\Gamma_p \equiv \mathcal{N}^l_{M(c)} = \mathcal{N}^{l'}_{M(c')}$, and $\mathbf{v}^l_{h,c, i} = \mathbf{v}^{l'}_{h,c', i}$, $i=1,\cdots,\mathcal{N}^\Gamma_p$.

\subsubsection{Truth problem statement}

We now formulate the truth problem in terms of bubble and port functions. As such, we introduce system's $\mathcal{N}_h^{\bb}$-dimensional bubble, $\mathcal{N}_h^{\Gamma}$-dimensional port(-lifted), and $\mathcal{N}_h$-dimensional FE spaces, respectively, given by $\mathcal{V}_h^\bb \equiv \oplus_{c \in \sys}\mathcal{V}_{h,c}^\bb$, $\mathcal{V}_h^{\Gamma} \equiv \oplus_{c \in \sys}\mathcal{V}_{h,c}^\gamma$, and $\mathcal{V}_h \equiv \left(\mathcal{V}_{h}^\bb \oplus \mathcal{V}^\Gamma_h \right) \cap \mathcal{V}$, where $\mathcal{N}_h^{\bb} = \sum_{c \in \sys} \mathcal{N}_{M(c)}^{\bb}$, $\mathcal{N}_h^{\Gamma} = \sum_{p \in \pset} \mathcal{N}_{p}^{\Gamma}$, and $\mathcal{N}_h = \mathcal{N}_h^{\bb} + \mathcal{N}_{h}^\Gamma$. Note that the intersection with $\mathcal{V}$ enforces essential boundary conditions. We further define $(\widehat{x}_{\wc,q}, \wrho_{\wc,q} )_{q=1}^{Q_\wc}$ $\forall \wc \in \lib$ as the truth quadrature rule in the reference domain $\widehat{\Omega}_\wc$ of each archetype component $\wc \in \lib$.  We may now state the truth problem: given $\mu = (\mu_c)_{c \in \sys} \in \mathcal{D}$, find $\{ {u}_{h,c}^{\bb}(\mu) \in \mathcal{V}_{h,c}^\bb \}_{c \in \sys}$ and $\{ \{ {u}_{h,c}^p(\mu) \in \mathcal{V}_{h,c}^p \}_{p \in \pset_{M(c)}} \}_{c \in \sys}$ such that, for all $\{ v_{h,c}^{\bb} \in \mathcal{V}_{h,c}^\bb \}_{c \in \sys}$ and $\{ \{ {v}_{h,c}^p \in \mathcal{V}_{h,c}^p \}_{p \in \pset_{M(c)}} \}_{c \in \sys}$,
\begin{equation}\label{eq:truthref}
    \begin{aligned}
    R_h(u_h(\mu), v_h;\mu) &\equiv \sum_{c \in \sys} \sum_{q=1}^{Q_{M(c)}} \wrho_{M(c),q} \widehat{r}_{M(c)} \biggl( \Bigl[ u^\bb_{h,c}(\mu) + \sum_{p \in \pset_{M(c)}} u^{p}_{h,c}(\mu) \Bigr] \circ \mathcal{G}_{c}(\widehat{x}_{M(c),q}; \mu_c), \\
    &\hspace{11.5em}\Bigl[ v^\bb_{h,c} + \sum_{p \in \pset_{M(c)}} v^{p}_{h,c} \Bigr] \circ \mathcal{G}_{c}(\widehat{x}_{M(c),q}; \mu_c); \widehat{x}_{M(c),q},\mu_c \biggr) = 0,
\end{aligned}
\end{equation}
where the system-level solution $u_h(\mu) \in \mathcal{V}_h$ is given by $u_h(\mu) =  \sum_{c \in \sys} \left[ u^\bb_{h,c}(\mu) + \sum_{p \in \pset_{M(c)}} u^{p}_{h,c}(\mu) \right]$. Similarly to the exact problem in~\eqref{eq:exact}, we assume~\eqref{eq:truthref} is well-posed for all $\mu \in \mathcal{D}$. We then evaluate the truth output
\begin{equation}\label{eq:truthoutput}
F_h(u_h(\mu);\mu) \equiv \sum_{c \in \sys} \sum_{q=1}^{Q_{M(c)}} \wrho_{M(c),q} \: \widehat{f}_{M(c)} \Big( \Bigl[u^\bb_{h,c}(\mu) + \sum_{p \in \pset_{M(c)}} u^{p}_{h,c}(\mu) \Bigr] \circ \mathcal{G}_{c}(\widehat{x}_{M(c),q}; \mu_c); \widehat{x}_{M(c),q}, \mu_c \Big).
\end{equation}

In practice, \eqref{eq:truthref} is solved using Newton's method. Given the $n$-th Newton iterate $u_h^{(n)}$, the $n+1$-st iterate is given by $u_h^{(n+1)} = u_h^{(n)} - \delta u_h^{(n)}$, where $\delta u_h^{(n)} \in \mathcal{V}_h$ is the Newton update that satisfies
\begin{equation}\label{eq:newtonupdate}
\begin{aligned}
R_h'(u_h^{(n)}, \delta u_h^{(n)},v_h;\mu) = R_h(u_h^{(n)},v_h;\mu) \quad \forall v_h \in \mathcal{V}_h,
\end{aligned}
\end{equation}
where $R_h'(u_h^{(n)}, \delta u_h^{(n)},v_h;\mu)$ is the G\^{a}teaux derivative of $R_h(\cdot,v_h;\mu)$ at $u_h^{(n)}$ in the direction of $\delta u_h^{(n)}$. We may appeal to~\eqref{eq:truthref} to obtain, $\forall w_h,z_h,v_h \in \mathcal{V}_h$ and $\forall \mu \in \mathcal{D}$,
\begin{equation}\nonumber
    \begin{aligned}
    R_h'(w_h, z_h,v_h;\mu) = \sum_{c \in \sys} \sum_{q=1}^{Q_{M(c)}} \wrho_{M(c),q} \: \widehat{r}\:'_{M(c)} \Bigl( &\Bigl[ w^\bb_{h,c} + \sum_{p \in \pset_{M(c)}} w^{p}_{h,c} \Bigr] \circ \mathcal{G}_{c}(\cdot; \mu_c), \\
    &\hspace{0em} \Bigl[ z^\bb_{h,c} + \sum_{p \in \pset_{M(c)}} z^{p}_{h,c} \Bigr] \circ  \mathcal{G}_{c}(\cdot; \mu_c), \\
    &\hspace{0em} \Bigl[ v^\bb_{h,c} + \sum_{p \in \pset_{M(c)}} v^{p}_{h,c} \Bigr] \circ \mathcal{G}_{c}(\cdot; \mu_c);\widehat{x}_{M(c),q},\mu_c \Bigr),
\end{aligned}
\end{equation}
where $\widehat{r}\:'_{\wc}(w_h,z_h,v_h;\widehat{x}_{\wc,q},\mu_\wc)$, $\wc \in \lib$, is the G\^{a}teaux derivative of $\widehat{r}_{\wc}(\cdot,v_h;\widehat{x}_{\wc,q},\mu_c)$ at $w_h$ in the direction of~$z_h$. 

\section{Hyperreduced reduced basis element method}
\label{sec:nlrbe}
In this section, we present our HRBE method, which uses component-wise RB and hyperreduction (i) to provide an accurate approximation of the truth problem~\eqref{eq:truthref} at a substantially reduced cost and (ii) to provide topological and parametric flexibility to assemble an arbitrary system in the online phase.

\subsection{RB problem formulation}
We now introduce an RB approximation of the truth problem~\eqref{eq:truthref}. We first construct an RB space for the bubble space of each component. To this end, we assume that, for each $c \in \sys$, the bubble solution $u_{h,c}^\bb(\mu)$ associated with~\eqref{eq:truthref} for any $\mu \in \mathcal{D}$ can be well-approximated in an $N_{M(c)}^\bb \ll \mathcal{N}_{M(c)}^\bb$-dimensional linear space. We then introduce, for each archetype component $\wc \in \lib$, an $N_{\wc}^\bb$-dimensional space $\what{\mathcal{V}}_{\rb,\wc}^\bb \subset \widehat{\mathcal{V}}_{h,\wc}^\bb$ spanned by an RB $\{ \widehat{\xi}_{\wc,i}^\bb\}_{i=1}^{N_{\wc}^\bb}$. (We defer the discussion of the computational procedure to construct RBs to Section~\ref{sec:training}; for now, we assume RBs are given.) We further define the $N_{\wc}$-dimensional space $\widehat{\mathcal{V}}_{\rb,\wc} \equiv \widehat{\mathcal{V}}_{\rb,\wc}^\bb \oplus \widehat{\mathcal{V}}_{h,\wc}^\gamma$, where $N_{\wc} = N_{\wc}^\bb + \mathcal{N}_\wc^\gamma$. Analogously, for each instantiated component $c \in \sys$, we introduce RB spaces $\mathcal{V}_{\rb,c}^\bb \equiv \left\{v = \widehat{v} \circ \mathcal{G}_c^{-1}(\cdot; \mu_c) \Big| \: \widehat{v} \in \widehat{\mathcal{V}}_{\rb,M(c)}^\bb \right\} \subset \mathcal{V}_{h,c}^\bb$ and $\mathcal{V}_{\rb,c} \equiv \mathcal{V}_{\rb,c}^\bb \oplus \mathcal{V}_{h,c}^\gamma$, so that we can express any $v_{\rb,c} \in \mathcal{V}_{\rb,c}$ as $v_{\rb,c} = v_{\rb,c}^\bb + v_{h,c}^\gamma$, where $v_{\rb,c}^\bb \in \mathcal{V}_{\rb,c}^\bb$ and $v_{h,c}^\gamma \in \mathcal{V}_{h,c}^\gamma$. 

We next define the system-level (global) RB space. We first introduce the bubble space for the system as the direct sum of component RB spaces: i.e., $\mathcal{V}_{\rb}^\bb \equiv \oplus_{c \in \sys} \mathcal{V}_{\rb,c}^\bb$. We then augment the space with the global port basis and enforce essential boundary conditions to obtain $\mathcal{V}_{\rb} \equiv \left(\mathcal{V}_{\rb}^\bb \oplus \mathcal{V}^\Gamma_h \right) \cap \mathcal{V}$. The dimensions of $\mathcal{V}_\rb^\bb$ and $\mathcal{V}_{\rb}$ are $N_{\rb}^\bb \equiv \sum_{c \in \sys} N^\bb_{M(c)}$ and $N_{\rb} \equiv N_{\rb}^\bb + \mathcal{N}^\Gamma_h$, respectively. Since $N_\wc^\bb \ll \mathcal{N}_\wc^\bb \: \: \forall \wc \in \lib$, we have $N_{\rb} \ll \mathcal{N}_h$. We then appeal to Galerkin projection to obtain the RB problem: given $\mu = (\mu_c)_{c \in \sys} \in \mathcal{D}$, find $\{ u_{\rb,c}^\bb(\mu) \in \mathcal{V}_{\rb,c}^\bb \}_{c \in \sys}$ and $\{ \{ u_{h,c}^p(\mu) \in \mathcal{V}_{h,c}^p \}_{p \in \pset_{M(c)}} \}_{c \in \sys}$ such that, for all $\{ v_{\rb,c}^\bb \in \mathcal{V}_{\rb,c}^\bb \}_{c \in \sys}$ and $\{ \{ v^p_{h,c} \in \mathcal{V}_{h,c}^p \}_{p \in \pset_{M(c)}} \}_{c \in \sys}$,
\begin{equation}\label{eq:rbtruth}
\begin{aligned}
R_{\rb} (u_{\rb}(\mu), v_{\rb};\mu ) = \sum_{c \in \sys} \sum_{q=1}^{Q_{M(c)}} \wrho_{M(c),q} \: \widehat{r}_{M(c)} \biggl( &\Bigl[ u^\bb_{\rb,c}(\mu) + \sum_{p \in \pset_{M(c)}} u^{p}_{h,c}(\mu) \Bigr] \circ \mathcal{G}_{c}(\what{x}_{M(c),q}; \mu_c), \\
&\Bigl[ v^\bb_{\rb,c} + \sum_{p \in \pset_{M(c)}} v^{p}_{h,c} \Bigr] \circ \mathcal{G}_{c}(\what{x}_{M(c),q}; \mu_c); \widehat{x}_{M(c),q},\mu_c \biggr) = 0,
\end{aligned}
\end{equation}
where the system-level solution $u_{\rb}(\mu) \in \mathcal{V}_{\rb}$ is given by $u_{\rb}(\mu) = \sum_{c \in \sys} \left[ u^\bb_{\rb,c}(\mu) + \sum_{p \in \pset_{M(c)}} u^{p}_{h,c}(\mu) \right]$, and evaluate the output $F_{\rb}(u_{\rb}(\mu);\mu)$. Here, $R_\rb(w,v;\mu) = R_h(w,v;\mu)$ and $F_\rb(w;\mu) = F_h(w;\mu)$ $\forall w,v \in \mathcal{V}_\rb, \: \forall \mu \in \mathcal{D}$, and hence the forms are evaluated using the truth quadrature rule. We again assume the RB problem~\eqref{eq:rbtruth} is well-posed for all $\mu \in \mathcal{D}$.

\begin{remark}
  \label{rem:no_port_red}
  In this work, we do not consider port reduction~\cite{eftang2013port, eftang2014port,smetana2015new}. Hence, the number of DoF in the RB system is bounded from below by the number of port DoF in the truth system, which ultimately limits the dimensionality reduction achieved by the present formulation, especially for systems with many ports and/or large ports. While recognizing the limitation, we focus on developing component-wise hyperreduction for this non-port-reduced system in this work and leave port-reduction to future work. 
\end{remark}

\subsection{HRBE problem formulation}
\label{subsec:hyper}
The computational cost of solving~\eqref{eq:rbtruth}, which uses the truth quadrature rule, depends on the underlying truth FE discretization, rendering the method not online efficient. To remedy this issue, we appeal to \emph{hyperreduction} techniques. Hyperreduction approaches in the RB literature fall into two classes. The first class of methods approximates integrands first and then integrates them. These methods use a number of empirically-derived basis functions to approximate the nonlinear terms in the integrands through a sparse interpolation/regression scheme and then integrate the approximated integrands.  Methods in this class include the gappy proper orthogonal decomposition (POD) method~\cite{everson1995karhunen}, the empirical interpolation method (EIM)~\cite{barrault2004empirical,grepl2007efficient}, the discrete EIM~\cite{chaturantabut2010nonlinear}, the first-order EIM~\cite{nguyen2023efficient, nguyen2024}, and the Gauss--Newton approximation tensor method~\cite{carlberg2011efficient}. The second class of hyperreduction methods directly approximates the integrals in the residual and output forms using a set of empirically-derived sparse element sampling or quadrature rule. Methods in this class include the optimal cubature method~\cite{An_2008_Cubature}, the energy-conserving mesh sampling and weighting method~\cite{farhat2014dimensional, farhat2015structure}, the empirical cubature method~\cite{Hernandez_2017_Empirical_Cubature}, and the EQP~\cite{Patera_2017_EQP, yano2019lp}. In the present study, we build on the EQP and its ability to construct quantitative control of the hyperreduction error in the solution (instead of the residual) and extend this capability to the component-based context.

The EQP constructs a set of empirical and sparse reduced quadrature (RQ) points and weights that approximate the integrals in the residual and output forms to a prescribed accuracy. The RQ points are a sparse subset of the truth quadrature points $\{\widehat{x}_{\wc,q} \}_{q=1}^{Q_\wc}$, $\wc \in \lib$, with re-weighted quadrature weights. We introduce $(\th{x}_{\wc,q}^r, \th{\rho}_{\wc,q}^r )_{q=1}^{\wtilde{Q}^r_\wc} \subset (\widehat{x}_{\wc,q}, \wrho_{\wc,q} )_{q=1}^{Q_\wc}$ as the residual RQ rule for each archetype component $\wc \in \lib$, where $\wtilde{Q}^r_\wc \ll Q_\wc$. We similarly introduce output functional RQ rule, $(\th{x}_{\wc,q}^f, \th{\rho}_{\wc,q}^f)_{q=1}^{\wtilde{Q}^f_\wc}  \subset (\widehat{x}_{\wc,q}, \wrho_{\wc,q} )_{q=1}^{Q_\wc}$, $\wc \in \lib$, where $\wtilde{Q}^f_\wc \ll Q_\wc$. (We defer the discussion of construction of these RQ rules in the offline phase to Section~\ref{sec:training}; for now, we assume the rules are given.) Given the RQ rules, the HRBE problem is stated as follows: given $\mu = (\mu_c)_{c \in \sys} \in \mathcal{D}$, find $\{ \widetilde{u}_{\rb,c}^{\bb}(\mu) \in \mathcal{V}_{\rb,c}^\bb \}_{c \in \sys}$ and $\{ \{ \wtilde{u}_{h,c}^p(\mu) \in \mathcal{V}_{h,c}^p \}_{p \in \pset_{M(c)}} \}_{c \in \sys}$ such that, for all $\{ v_{\rb,c}^{\bb} \in \mathcal{V}_{\rb,c}^\bb \}_{c \in \sys}$ and $\{ \{ {v}_{h,c}^p \in \mathcal{V}_{h,c}^p \}_{p \in \pset_{M(c)}} \}_{c \in \sys}$,
\begin{equation}\label{eq:nlrbepr}
\begin{aligned}
\widetilde{R}_{\rb} (\widetilde{u}_{\rb}(\mu), {v}_{\rb};\mu ) \equiv \sum_{c \in \sys} \sum_{q=1}^{\wtilde{Q}_{M(c)}^r} \th{\rho}_{M(c),q}^r \: \widehat{r}_{M(c)} \biggl( &\Bigl[ \widetilde{u}^{\bb}_{\rb,c}(\mu) + \sum_{p \in \pset_{M(c)}} \wtilde{u}^{p}_{h,c}(\mu) \Bigr] \circ {\mathcal{G}}_{c}(\th{x}_{M(c),q}^r; \mu_c),  \\
&\Bigl[{v}^{\bb}_{\rb,c} + \sum_{p \in \pset_{M(c)}} v^{p}_{h,c} \Bigr] \circ {\mathcal{G}}_{c}(\th{x}_{M(c),q}^r; \mu_c); \th{x}_{M(c),q}^r,\mu_c \biggr) = 0,
\end{aligned}
\end{equation}
and evaluate the approximate output
\begin{equation}\label{eq:outputnlrbe}
\begin{aligned}
\widetilde{F}_{\rb}(\widetilde{u}_{\rb} (\mu);\mu) \equiv \sum_{c \in \sys} \sum_{q=1}^{\wtilde{Q}_{M(c)}^f} \th{\rho}_{M(c),q}^f \: \widehat{f}_{M(c)} \biggl( \Bigl[ \widetilde{u}^{\bb}_{\rb,c}(\mu) + \sum_{p \in \pset_{M(c)}} \wtilde{u}^{p}_{h,c}(\mu) \Bigr] \circ {\mathcal{G}}_{c}(\th{x}_{M(c),q}^f; \mu_c); \th{x}_{M(c),q}^f, \mu_c \biggr).
\end{aligned}
\end{equation}
Owing to $N_{\rb} \ll \mathcal{N}_h$, $\wtilde{Q}^r_c \ll Q_{M(c)}$, and $\wtilde{Q}_c^f \ll Q_{M(c)}$ $\forall c \in \sys$, solving the hyperreduced RB problem~\eqref{eq:nlrbepr} and approximating the output~\eqref{eq:outputnlrbe} can be carried out significantly more efficiently than their corresponding counterparts in the truth problem,~\eqref{eq:truthref} and~\eqref{eq:truthoutput}, respectively. Sufficient conditions for the well-posedness of the hyperreduced RB problem~\eqref{eq:outputnlrbe} will be provided in Proposition~\ref{prop:suffcond}.

\begin{remark}
  For each archetype component, the bubble RB and RQ rules are calculated and stored in the library \emph{a priori} in the offline phase. Therefore, in the online phase, once we determine the connectivity of instantiated components and form a system, we can rapidly assemble the system's reduced residual and Jacobian and solve the HRBE system~\eqref{eq:nlrbepr} without RB and RQ retraining. In other words, the HRBE system results from assembling hyperreduced components trained in the offline phase, and \emph{not} from applying hyperreduction to an online-assembled RB system, which could not be performed in an online efficient manner.
\end{remark}

\section{Component-wise RB and RQ training}
\label{sec:training}
The two primary ingredients of the HRBE method presented in Section~\ref{sec:nlrbe} are the RB  $\{ \widehat{\xi}_{\wc,i}^\bb\}_{i=1}^{N_{\wc}^\bb}$ of the bubble spaces $\widehat{\mathcal{V}}_{\rb,\wc}^\bb$ and the RQ rules $(\th{x}_{\wc,q}^r, \th{\rho}_{\wc,q}^r )_{q=1}^{\wtilde{Q}^r_\wc}$ and $(\th{x}_{\wc,q}^f, \th{\rho}_{\wc,q}^f )_{q=1}^{\wtilde{Q}^f_\wc}$ for all archetype components $\wc \in \lib$. In this section, we outline the procedures to construct these essential elements. 

\subsection{Generation of archetype-component training solutions}
\label{subsec:train_set_gen}
For each archetype component, we use an empirical training procedure to deduce the shape and magnitude of its anticipated solution and boundary conditions. To this end, we introduce for each archetype component $\wc \in \lib$, a parameter training set $\Xi^{\train}_\wc \equiv \{ \mu^{\train}_{\wc,n} \in \widehat{\mathcal{D}}_\wc \}_{n=1}^{N^{\train}_\wc}$ with a size of $N^{\train}_\wc$. For this archetype component, we compose $N_{\sample}$ sample \emph{subsystems} by connecting it through its $n^\gamma_\wc$ local ports to other randomly selected components from the library. This component is connected to other components through its local ports with a probability of $\beta$. We then assign random parameter values to each component in the assembled subsystems from their respective parameter training sets. We next apply random independent constant Dirichlet boundary conditions, with uniform density, to all nonshared global ports. We finally solve the truth problem for each subsystem, extract the truth solutions on the target component to form a state snapshot set $U_{h,\wc}^{\train} \equiv \{ u_{h,\wc,n}^\train \}_{n=1}^{N_\sample}$ designated for this component. The fundamental assumption underpinning this process is that the generated set of snapshot solutions sufficiently represents the set of all potential solutions and boundary conditions the component may experience in an actual system configuration. Algorithm~\ref{alg:emptr} provides an outline of the empirical process to generate snapshot solutions for archetype components.

\begin{algorithm}[!tb]
\caption{Generating training data for RB and RQ construction of archetype components.}\label{alg:emptr}
\SetAlgoLined
\KwInput{Number of sample subsystems $N_{\sample}$; probability of port connection $0 \leq \beta \leq 1$}
\KwOutput{The set of snapshot solutions $U_{h,\wc}^{\train}$ $\forall \wc \in \lib$}
    \For{$\wc \in \lib$}{
    	$U_{h,\wc}^{\train} = \emptyset$\;
      \For{$n=1, \cdots, N_{\sample}$}{
        \tcp{Assemble subsystem $\mathcal{C}_{\rm sub}$ and extract the solution associated with $\wc$}
    		\For{$\wp \in \pset_\wc$}{
    			Connect the archetype component $\wc$ through its $\wp$-th port to another component in the library with a probability of $\beta$\; 
    		}
    		Assign parameter value $\mu_c$ drawn uniform-randomly from $\mathcal{D}_c$ to each component $c \in \mathcal{C}_{\rm sub}$\; 
    		Assign uniform-random constant Dirichlet boundary conditions to each nonshared global port\; 
    		Solve the truth problem for the composed subsystem $\mathcal{C}_{\rm sub}$\; \label{ln:solve}
    		Extract the solution $u_{h,\wc,n}^{\train}$ on component $\wc$\;
    		$U_{h,\wc}^{\train} \leftarrow U_{h,\wc}^{\train} \cup u_{h,\wc,n}^{\train}$\;
    	}
    }
\end{algorithm}

\subsection{Component-wise RB construction}
\label{subsec:rbptraining}

For each archetype component $\wc \in \lib$, we decompose snapshot solutions in the training set $U_{h,\wc}^\train$ into their bubble and port solutions as in~\eqref{eq:compsol}. The bubble solutions are added to a training set $U^{\train,\bb}_{h,\wc}$ considered for this component. We then apply the POD to construct an RB $\{ \widehat{\xi}_{\wc,i}^\bb\}_{i=1}^{N_{\wc}^\bb}$ for the bubble space $\widehat{\mathcal{V}}_{\rb,\wc}^\bb$.
 
\subsection{Component-wise hyperreduction: BRR theory}
\label{subsec:eqptraining_brr}
We now present an extension of the original EQP~\cite{Patera_2017_EQP, yano2019lp} to the component-based context. To this end, we first introduce the BRR theorem~\cite{caloz1997numerical} specialized for the Euclidean space. 

\begin{lemma}[Brezzi–Rappaz–Raviart theorem]\label{lemma:brr}
    Given an $N$-dimensional Euclidean space $\mathbb{R}^N$, we introduce a $C^1$ mapping $G: \mathbb{R}^N \rightarrow \mathbb{R}^N$, $\mathbf{v} \in \mathbb{R}^N$ such that the Jacobian $DG(\mathbf{v}) \in \mathbb{R}^{N \times N}$ is nonsingular, and constants $\varepsilon$, $\delta$, and $L(\alpha)$ such that
    \begin{align}
    \left\|G(\mathbf{v}) \right\|_2 &\leq \varepsilon, \nonumber\\
    \left\| DG^{-1}(\mathbf{v}) \right\|_2 &\leq \delta, \nonumber\\
    \sup_{\mathbf{w} \in \bar{B}(\mathbf{v}, \alpha)}  \left\| DG(\mathbf{v}) -  DG(\mathbf{w}) \right\|_2 &\leq L(\alpha),\nonumber
    \end{align}
    where $\bar{B}(\mathbf{v}, \alpha) \equiv \{ \mathbf{z}: \left\| \mathbf{z} - \mathbf{v} \right\|_2 \leq \alpha \}$. Assume $2 \delta L(2 \delta \varepsilon) \leq 1$. Then, for all $\lambda \geq 2 \delta \varepsilon$ such that $\delta L(\lambda) < 1$, there exists a unique $\mathbf{u} \in \mathbb{R}^N$ that satisfies $G(\mathbf{u}) = 0$ in the ball $\bar{B}(\mathbf{v}, 2 \delta \varepsilon)$ and $DG(\mathbf{u}) \in \mathbb{R}^{N \times N}$ is invertible and satisfies
    \begin{equation}\label{eq:brrjac}
    \left\| DG^{-1}(\mathbf{u}) \right\|_2 \leq 2 \left\| DG^{-1}(\mathbf{v}) \right\|_2 \leq 2 \delta.
    \end{equation}
    Additionally, 
    \begin{equation}\label{eq:brrerr}
    \left\| \mathbf{w} - \mathbf{u} \right\|_2 \leq 2 \left\| DG^{-1}(\mathbf{v}) \right\|_2 \left\| G(\mathbf{w}) \right\|_2  \leq 2 \delta \left\| G(\mathbf{w}) \right\|_2 \quad \forall \mathbf{w} \in \bar{B}(\mathbf{v}, 2 \delta \varepsilon).
    \end{equation}
    \end{lemma} 
\begin{proof}
    See~\cite{caloz1997numerical}.
 \end{proof}

\begin{corollary}[Effectivity bound]\label{coro:effectivity}
  For all $\mathbf{w} \in \bar{B}(\mathbf{v}, 2 \delta \varepsilon)$, the effectivity of the error bound $\Delta(\mathbf{w}) \equiv 2 \delta \left\| G(\mathbf{w}) \right\|_2$ is bounded by
    \begin{equation}\label{eq:effectivity}
      \eta(\mathbf{w})
      \equiv \frac{\Delta(\mathbf{w})}{\left\| \mathbf{w} - \mathbf{u} \right\|_2}
      \leq 2 \delta \Big( \bar{L}(\mathbf{w}) +  \| DG(\mathbf{u}) \|_2 \Big),
    \end{equation}
    where 
    \begin{equation}\label{eq:Lbar}
        \bar{L}({\mathbf{z}}) \equiv \sup_{{\mathbf{z}} \in \bar{B}(\mathbf{w}, \| \mathbf{u} - \mathbf{w} \|_2)}  \left\| DG(\mathbf{u}) -  DG({\mathbf{z}}) \right\|_2.
    \end{equation}
 \end{corollary}
 \begin{proof}
    We first consider the Taylor expansion of $G(\cdot)$ about $\mathbf{w} \in \mathbb{R}^N$,
    \begin{equation}
        \begin{aligned}
        G(\mathbf{u}) &= G(\mathbf{w}) + \int_{0}^1 DG (\mathbf{w} + t(\mathbf{u} - \mathbf{w}) ) (\mathbf{u} - \mathbf{w}) dt \\
        &= G(\mathbf{w}) + \int_{0}^1 \Big[DG (\mathbf{w} + t(\mathbf{u} - \mathbf{w}) ) - DG(\mathbf{u}) \Big] (\mathbf{u} - \mathbf{w}) dt + DG(\mathbf{u}) (\mathbf{u} - \mathbf{w}).
        \end{aligned}
    \end{equation}
    We next appeal to $G(\mathbf{u}) = 0$ and \eqref{eq:Lbar} to obtain $\| G(\mathbf{w}) \|_2 \leq \Big( \bar{L}(\mathbf{w}) +  \| DG(\mathbf{u}) \|_2 \Big) \: \| \mathbf{u} - \mathbf{w} \|_2$. We then multiply both sides by $2\delta/\| \mathbf{u} - \mathbf{w} \|_2$ to obtain~\eqref{eq:effectivity}.
 \end{proof}

\begin{remark}
  \label{rem:brr_applicability}
  Mathematically, the BRR theorem holds for any nonlinear mapping that satisfies the conditions in Lemma~\ref{lemma:brr}. In the context of numerical approximation of parametrized PDEs, the theorem is applicable as long as the approximate solution is sufficiently close to the solution and the Jacobian is nonsingular and Lipschitz continuous. In practice, the BRR theorem can be (and has been) applied to PDEs including Burgers' equation~\cite{VEROY2003619, Yano2014Burgers}, (low-Reynolds-number) Navier--Stokes equations~\cite{veroy2005certified, Yano2014Boussinesq}, hyperelasticity~\cite{yano2019lp}, and heat transfer (as considered in Section~\ref{sec:res}). However, for certain problems, such as convection-dominated high-Reynolds-number Navier--Stokes equations, the BRR error bound can be overly conservative and may offer limited practical utility.
 \end{remark}

 We now specialize this lemma to the component-based context to facilitate the development of our component-wise hyperreduction scheme. As such, for each instantiated component $c \in \sys$, we introduce the algebraic RB residual and Jacobian ${\mathbf{R}}_{\rb,c}: \mathbb{R}^{N_{M(c)}} \times \mathcal{D}_c \to \mathbb{R}^{N_{M(c)}}$ and ${\mathbf{J}}_{\rb,c}: \mathbb{R}^{N_{M(c)}} \times \mathcal{D}_c \to \mathbb{R}^{N_{M(c)} \times N_{M(c)}}$ as well as the algebraic hyperreduced RB residual and Jacobian $\wtilde{\mathbf{R}}_{\rb,c}: \mathbb{R}^{N_{M(c)}} \times \mathcal{D}_c \to \mathbb{R}^{N_{M(c)}}$ and $\wtilde{\mathbf{J}}_{\rb,c}: \mathbb{R}^{N_{M(c)}} \times \mathcal{D}_c \to \mathbb{R}^{N_{M(c)} \times N_{M(c)}}$. For conciseness, we omit the explicit expressions of these quantities here and provide them instead in Appendix~\ref{app:resandjacrb}. Additionally, we introduce $\mathbf{P}_{\rb,c}: \mathbb{R}^{{N}_{M(c)}} \rightarrow \mathbb{R}^{{N}_\rb}$ as the linear extension operator that maps the components' RB DoF to the assembled system's DoF. 
 
 At the system level, we introduce $\mathbf{R}_{\rb}: \mathbb{R}^{N_{\rb}} \times \mathcal{D} \rightarrow \mathbb{R}^{N_{\rb}}$ and $\mathbf{J}_{\rb}: \mathbb{R}^{N_{\rb}} \times \mathcal{D} \rightarrow \mathbb{R}^{N_{\rb} \times {N_{\rb}}}$ as, respectively, the algebraic residual and Jacobian of the (truth-quadrature) RB problem~\eqref{eq:rbtruth}. We denote the generalized coordinates of the solution by $\mathbf{u}_{\rb}(\mu) \in \mathbb{R}^{N_{\rb}}$. We similarly denote the algebraic residual and Jacobian of the hyperreduced RB problem~\eqref{eq:nlrbepr} by $\widetilde{\mathbf{R}}_{\rb}: \mathbb{R}^{N_{\rb}} \times \mathcal{D} \rightarrow \mathbb{R}^{N_{\rb}}$ and $\widetilde{\mathbf{J}}_{\rb}: \mathbb{R}^{N_{\rb}} \times \mathcal{D} \rightarrow \mathbb{R}^{N_{\rb} \times N_{\rb}}$, respectively, and denote the generalized coordinates of the solution by $\widetilde{\mathbf{u}}_{\rb}(\mu) \in \mathbb{R}^{N_{\rb}}$. The system-level residual and Jacobian $\forall \mathbf{w}_\rb \in \mathbb{R}^{N_\rb}$ and $\forall \mu \in \mathcal{D}$ can be obtained through
 \begin{align}
    {\mathbf{R}}_{\rb} (\mathbf{w}_{\rb}; \mu) = \sum_{c \in \sys} \mathbf{P}_{\rb,c} {\mathbf{R}}_{\rb,c}(\mathbf{P}_{\rb,c}^T \mathbf{w}_{\rb}; \mu_c), \quad {\mathbf{J}}_{\rb} (\mathbf{w}_{\rb}; \mu)  = \sum_{c \in \sys} \mathbf{P}_{\rb,c} {\mathbf{J}}_{\rb,c}(\mathbf{P}_{\rb,c}^T \mathbf{w}_{\rb}; \mu_c) \mathbf{P}_{\rb,c}^T, \nonumber \\
 \widetilde{\mathbf{R}}_{\rb} (\mathbf{w}_{\rb}; \mu) = \sum_{c \in \sys} \mathbf{P}_{\rb,c} \widetilde{\mathbf{R}}_{\rb,c}(\mathbf{P}_{\rb,c}^T \mathbf{w}_{\rb}; \mu_c), \quad \widetilde{\mathbf{J}}_{\rb} (\mathbf{w}_{\rb}; \mu)  = \sum_{c \in \sys} \mathbf{P}_{\rb,c} \widetilde{\mathbf{J}}_{\rb,c}(\mathbf{P}_{\rb,c}^T \mathbf{w}_{\rb}; \mu_c) \mathbf{P}_{\rb,c}^T. \nonumber
 \end{align}
 
 \begin{proposition}\label{prop:suffcond}
    For a system $\sys$ and given $\mu \in \mathcal{D}$, we introduce $\bar{u}_\rb(\mu) \in \mathcal{V}_\rb$ and its associated generalized coordinates $\bar{\mathbf{u}}_{\rb}(\mu) \in \mathbb{R}^{N_{\rb}}$ such that 
    \begin{equation}\label{eq:cond}
    \left\| \mathbf{u}_{\rb}(\mu) -  \bar{\mathbf{u}}_{\rb}(\mu) \right\|_2 \leq \bar{\varepsilon}
    \end{equation}
    for an $\bar{\varepsilon} \geq 0$ and $\mathbf{J}_{\rb}(\bar{\mathbf{u}}_{\rb}(\mu);\mu)$ is nonsingular. We further introduce $\sigma \equiv \sigma_{\mn} (\mathbf{J}_{\rb}\left(\bar{\mathbf{u}}_{\rb}(\mu);\mu) \right)$, where $\sigma_{\mn}(\cdot)$ denotes the minimum singular value of its argument. We suppose for some $\delta_{R_c} \geq 0$ and $\delta_{J_c} \geq 0$, $c \in \sys$, such that $\sum_{c \in \sys} N_{M(c)} \delta_{J_c} < \sigma$, the following inequalities hold:
    \begin{align}
    \left\| \widetilde{\mathbf{R}}_{\rb,c}(\mathbf{P}^T_{\rb,c} \bar{\mathbf{u}}_{\rb}(\mu);\mu_c) \right\|_\infty &\leq \delta_{R_c} &&\forall c \in \sys, \label{eq:const1}\\
    \left\| {\mathbf{J}}_{\rb,c}(\mathbf{P}^T_{\rb,c} \bar{\mathbf{u}}_{\rb}(\mu);\mu_c) - \widetilde{\mathbf{J}}_{\rb,c}(\mathbf{P}^T_{\rb,c} \bar{\mathbf{u}}_{\rb,c}(\mu);\mu_c) \right\|_\mx &\leq \delta_{J_c} &&\forall c \in \sys \label{eq:const2}, 
    \end{align}
    where for any $c \in \sys$, $\| \mathbf{A} \|_\mx \equiv \max_{i,j \in \{ 1,\cdots,N_{M(c)} \}}$  $|A_{i,j}|$ for $\mathbf{A} \in \mathbb{R}^{N_{M(c)} \times N_{M(c)}}$. We also introduce
    \begin{equation}\label{eq:lalpha}
    L(\alpha) \equiv 2 \sup_{\mathbf{w} \in \bar{B}(\bar{\mathbf{u}}_{\rb}(\mu), \alpha)}  \left\| \mathbf{J}_{\rb}^{-1}(\bar{\mathbf{u}}_{\rb}(\mu);\mu) \widetilde{\mathbf{J}}_{\rb}(\mathbf{w};\mu) - \mathbf{I} \right\|_2
    \end{equation}
    and assume 
    \begin{equation}\label{eq:lcond}
    L( \bar{\alpha}) \leq \frac{\sigma - \sum_{c \in \sys} N_{M(c)} \delta_{J_c}}{2 \sigma},
    \end{equation}
    where $\bar{\alpha} = 2  {\sum_{c \in \sys} \sqrt{N_{M(c)}} \delta_{R_c} } / (\sigma - \sum_{c \in \sys} N_{M(c)} \delta_{J_c})$. Then, for all $\lambda \geq \bar{\alpha}$, there exists a unique solution $\widetilde{\mathbf{u}}_{\rb}(\mu) \in \mathbb{R}^{N_{\rb}}$ such that $\widetilde{\mathbf{R}}_{\rb}(\widetilde{\mathbf{u}}_{\rb}; \mu) = 0$ in the ball $\bar{B}(\bar{\mathbf{u}}_{\rb}(\mu), \lambda)$, where $L(\lambda) \leq (\sigma - \sum_{c \in \sys} N_{M(c)} \delta_{J_c}) / \sigma$. Furthermore,
    \begin{equation}\label{eq:properrabs}
    \begin{aligned}
    \lr{ \mathbf{u}_{\rb}(\mu) - \widetilde{\mathbf{u}}_{\rb}(\mu) }_2 &\leq \bar{\alpha} +\bar{\varepsilon}.
    \end{aligned}
    \end{equation}
\end{proposition}
\begin{proof}
    For notational brevity, we suppress $\mu$ and $\mu_c \: \forall c \in \sys$ throughout the proof. Referring to Lemma~\ref{lemma:brr}, we set $G(\cdot) \equiv \mathbf{J}_{\rb}^{-1}(\bar{\mathbf{u}}_{\rb}) \widetilde{\mathbf{R}}_{\rb}(\cdot)$ and $\mathbf{v} \equiv \bar{\mathbf{u}}_{\rb}$. We observe that
    \begin{equation}\nonumber
    \begin{aligned}
    \lr{G(\mathbf{v})}_2 &= \lr{ \mathbf{J}_{\rb}^{-1}(\bar{\mathbf{u}}_{\rb}) \widetilde{\mathbf{R}}_{\rb}(\bar{\mathbf{u}}_{\rb}) }_2 \leq \lr{ \mathbf{J}_{\rb}^{-1}(\bar{\mathbf{u}}_{\rb}) }_2 \lr{ \widetilde{\mathbf{R}}_{\rb}(\bar{\mathbf{u}}_{\rb}) }_2 \leq \frac{{ \sum_{c \in \sys} \lr{ \widetilde{\mathbf{R}}_{\rb,c}(\mathbf{P}^T_{\rb,c} \bar{\mathbf{u}}_{\rb}) }_2}} {\sigma} \\
    &\leq \frac{{ \sum_{c \in \sys} \sqrt{N_{M(c)}} \lr{ \widetilde{\mathbf{R}}_{\rb,c}(\mathbf{P}^T_{\rb,c} \bar{\mathbf{u}}_{\rb}) }_\infty}} {\sigma} \leq \frac{{\sum_{c \in \sys} \sqrt{N_{M(c)}} \delta_{R_{M(c)}}}}{\sigma},
    \end{aligned}
    \end{equation}
    where the second inequality follows from the component-wise decomposition of the residual and the matrix norm relation $\| \mathbf{J}_{\rb}^{-1}(\bar{\mathbf{u}}_\rb) \|_2 = \sigma_{\mn}^{-1}(\mathbf{J}_{\rb}(\bar{\mathbf{u}}_\rb)) = \sigma^{-1}$, the third inequality follows from the relationship between $\| \cdot \|_2$ and $\| \cdot \|_\infty$, and the last inequality follows from condition~\eqref{eq:const1}.
    Hence, we set $\varepsilon \equiv {\sum_{c \in \sys} \sqrt{N_{M(c)}} \delta_{R_{M(c)}}} / \sigma$ in Lemma~\ref{lemma:brr}. Moreover, we have
    \begin{equation}\label{eq:brr_inter}
      \begin{aligned}
        \lr{ \mathbf{I} - \mathbf{J}_\rb^{-1} (\bar{\mathbf{u}}_\rb) \widetilde{\mathbf{J}}_\rb (\bar{\mathbf{u}}_\rb) }_2 &\leq \lr{\mathbf{J}_\rb^{-1} (\bar{\mathbf{u}}_\rb)}_2 \lr{ {\mathbf{J}}_{\rb}(\bar{\mathbf{u}}_{\rb}) - \widetilde{\mathbf{J}}_{\rb}(\bar{\mathbf{u}}_{\rb}) }_2 
        = \frac{1}{\sigma} \lr{ {\mathbf{J}}_{\rb}(\bar{\mathbf{u}}_{\rb}) - \widetilde{\mathbf{J}}_{\rb}(\bar{\mathbf{u}}_{\rb}) }_2  \\
          &= \frac{1}{\sigma}\lr{\sum_{c \in \sys} \mathbf{P}_{\rb,c} \left( {\mathbf{J}}_{\rb,c}(\mathbf{P}_{\rb,c}^T \bar{\mathbf{u}}_{\rb}) - \widetilde{\mathbf{J}}_{\rb,c}(\mathbf{P}_{\rb,c}^T \bar{\mathbf{u}}_{\rb}) \right) \mathbf{P}_{\rb,c}^T}_2  \\
    &\leq \frac{1}{\sigma} \sum_{c \in \sys} \lr{ {\mathbf{J}}_{\rb,c}(\mathbf{P}_{\rb,c}^T \bar{\mathbf{u}}_{\rb}) - \widetilde{\mathbf{J}}_{\rb,c}(\mathbf{P}_{\rb,c}^T \bar{\mathbf{u}}_{\rb}) }_2\\
          &\leq \frac{1}{\sigma}\sum_{c \in \sys} N_{M(c)} \lr{ {\mathbf{J}}_{\rb}(\mathbf{P}_{\rb,c}^T \bar{\mathbf{u}}_{\rb}) - \widetilde{\mathbf{J}}_{\rb}(\mathbf{P}_{\rb,c}^T \bar{\mathbf{u}}_{\rb}) }_\mx
          \leq \frac{1}{\sigma} \sum_{c \in \sys} N_{M(c)} \delta_{J_{M(c)}} < 1, 
    \end{aligned}
    \end{equation}
    where the first equality follows from the definition $\| \mathbf{J}_{\rb}^{-1}(\bar{\mathbf{u}}_\rb) \|_2 = \sigma_{\mn}^{-1}(\mathbf{J}_{\rb}(\bar{\mathbf{u}}_\rb)) = \sigma^{-1}$, the second equality follows from the component-wise decomposition of the Jacobian, the second inequality follows from the triangle inequality, the third inequality follows from the relation $\| \mathbf{A} \|_2 \leq N_\rb \| \mathbf{A} \|_\mx$ $\forall \mathbf{A} \in \mathbb{R}^{N_\rb \times N_\rb}$,  the fourth inequality follows from condition~\eqref{eq:const2}, and the last inequality follows from the assumption $\sum_{c \in \sys} N_{M(c)} \delta_{J_c} < \sigma$. It hence follows that
    \begin{equation}\nonumber
        \begin{aligned}
        \lr{ DG^{-1} (\mathbf{v}) }_2 &= \lr{ ( \mathbf{J}_\rb^{-1} (\bar{\mathbf{u}}_\rb) \widetilde{\mathbf{J}}_\rb (\bar{\mathbf{u}}_\rb))^{-1} }_2 =  \lr{ ( \mathbf{I} + \mathbf{J}_\rb^{-1} (\bar{\mathbf{u}}_\rb) \widetilde{\mathbf{J}}_\rb (\bar{\mathbf{u}}_\rb) - \mathbf{I})^{-1} }_2 \\
        &\leq \frac{1}{1 - \lr{ \mathbf{J}_\rb^{-1} (\bar{\mathbf{u}}_\rb) \widetilde{\mathbf{J}}_\rb (\bar{\mathbf{u}}_\rb) - \mathbf{I} }_2 } \leq \frac{\sigma}{\sigma - \sum_{c \in \sys} N_{M(c)} \delta_{J_c}},
        \end{aligned}
    \end{equation}
    where the first inequality follows from the Banach lemma which states $\forall \mathbf{A} \in \mathbb{R}^{N_\rb \times N_\rb}$ with $\| \mathbf{A} \|_2 < 1$, $(\mathbf{I} + \mathbf{A})^{-1}$ exists and satisfies $\| (\mathbf{I} + \mathbf{A})^{-1} \|_2 \leq (1 - \| \mathbf{A} \|_2)^{-1}$, and the last inequality follows from~\eqref{eq:brr_inter}. 
    We hence set $\delta \equiv \sigma / (\sigma - \sum_{c \in \sys} N_{M(c)} \delta_{J_c})$ in Lemma~\ref{lemma:brr}. We in addition note that 
    \begin{equation}\nonumber
    \begin{aligned}
    \sup_{\mathbf{w} \in \bar{B}(\mathbf{v}, \alpha)}  \lr{ DG(\mathbf{v}) -  DG(\mathbf{w}) }_2 &= \sup_{\mathbf{w} \in \bar{B}(\bar{\mathbf{u}}_\rb, \alpha)} \lr{ \mathbf{J}_\rb^{-1} (\bar{\mathbf{u}}_\rb) \widetilde{\mathbf{J}}_\rb (\bar{\mathbf{u}}_\rb) - \mathbf{J}_\rb^{-1} (\bar{\mathbf{u}}_\rb) \widetilde{\mathbf{J}}_\rb ( \mathbf{w} ) }_2 \\
    &= \sup_{\mathbf{w} \in \bar{B}(\bar{\mathbf{u}}_\rb, \alpha)} \lr{ \mathbf{J}_\rb^{-1} (\bar{\mathbf{u}}_\rb) \widetilde{\mathbf{J}}_\rb (\bar{\mathbf{u}}_\rb) - \mathbf{I} + \mathbf{I} - \mathbf{J}_\rb^{-1} (\bar{\mathbf{u}}_\rb) \widetilde{\mathbf{J}}_\rb ( \mathbf{w} ) }_2 \\
    &\leq \lr{ \mathbf{J}_\rb^{-1} (\bar{\mathbf{u}}_\rb) \widetilde{\mathbf{J}}_\rb (\bar{\mathbf{u}}_\rb) - \mathbf{I} }_2 + \sup_{\mathbf{w} \in \bar{B}(\bar{\mathbf{u}}_\rb, \alpha)} \lr{ \mathbf{J}_\rb^{-1} (\bar{\mathbf{u}}_\rb) \widetilde{\mathbf{J}}_\rb ( \mathbf{w} ) - \mathbf{I} }_2 \\
    &\leq 2 \sup_{\mathbf{w} \in \bar{B}(\bar{\mathbf{u}}_\rb, \alpha)} \lr{ \mathbf{J}_\rb^{-1} (\bar{\mathbf{u}}_\rb) \widetilde{\mathbf{J}}_\rb ( \mathbf{w} ) - \mathbf{I} }_2,
    \end{aligned}
    \end{equation}
    where the first inequality follows from the triangle inequality, and the last inequality follows from $\bar{\mathbf{u}}_\rb \in \bar{B}( \bar{\mathbf{u}}_\rb, \alpha)$. We hence set $L(\alpha) \equiv 2 \sup_{\mathbf{w} \in \bar{B}(\bar{\mathbf{u}}_\rb, \alpha)} \lr{ \mathbf{J}_\rb^{-1} (\bar{\mathbf{u}}_\rb) \widetilde{\mathbf{J}}_\rb ( \mathbf{w} ) - \mathbf{I} }_2$ in \eqref{eq:lalpha} in Lemma~\ref{lemma:brr}.

    Having defined $\epsilon$, $\delta$, and $L(\alpha)$ in the BRR theorem in Lemma~\ref{lemma:brr} for the HRBE method, we now apply the BRR theorem. If $2 \delta L (2 \delta \varepsilon) = 2  \sigma L(\bar{\alpha}) / ( \sigma - \sum_{c \in \sys} N_{M(c)} \delta_{J_c}) \leq 1$, we readily deduce, for all $\lambda \geq  2 \delta \varepsilon = \bar{\alpha}$ such that $L(\lambda) < 1 / \delta = (\sigma - \sum_{c \in \sys} N_{M(c)} \delta_{J_c}) / \sigma $, the existence of a unique solution $\mathbf{z} \in \mathbb{R}^{N_\rb}$ that satisfies $ G(\mathbf{z}) = \mathbf{J}_{\rb}^{-1}(\bar{\mathbf{u}}_{\rb}) \widetilde{\mathbf{R}}_{\rb}(\mathbf{z}) = 0$ in the ball $\bar{B}(\bar{\mathbf{u}}_\rb, \lambda )$. Since $\widetilde{\mathbf{u}}_\rb$ satisfies $\widetilde{\mathbf{R}}_\rb (\widetilde{\mathbf{u}}_\rb) = 0$, we conclude it is indeed the unique solution to both $G(\cdot) = 0$ and $\widetilde{\mathbf{R}}_\rb (\cdot) = 0$. Moreover, we set $\mathbf{w} \equiv \widetilde{\mathbf{u}}_\rb = \mathbf{v}$ in \eqref{eq:brrerr} to obtain
    \begin{equation}\nonumber
        \begin{aligned}
        \lr{ \mathbf{u}_\rb - \widetilde{\mathbf{u}}_\rb }_2 &\leq \lr{ \mathbf{u}_\rb - \bar{\mathbf{u}}_\rb }_2 + \lr{\bar{\mathbf{u}}_\rb -  \widetilde{\mathbf{u}}_\rb}_2 \leq\bar{\varepsilon} + 2 \delta \lr{ \mathbf{J}_{\rb}^{-1}(\bar{\mathbf{u}}_{\rb}) \widetilde{\mathbf{R}}_{\rb}(\bar{\mathbf{u}}_{\rb})  }_2 \\
        &\leq\bar{\varepsilon} + \frac{2 \sigma}{\sigma - \sum_{c \in \sys} N_{M(c)} \delta_{J_c}} \lr{\mathbf{J}_{\rb}^{-1}(\bar{\mathbf{u}}_{\rb})  }_2 \lr{ \widetilde{\mathbf{R}}_{\rb}(\bar{\mathbf{u}}_{\rb}) }_2 \\
        &\leq \bar{\varepsilon} + 2\frac{{\sum_{c \in \sys} \sqrt{N_{M(c)}} \delta_{R_{M(c)}}}}{\sigma - \sum_{c \in \sys} N_{M(c)} \delta_{J_c}} =\bar{\varepsilon} + \bar{\alpha},
        \end{aligned}
    \end{equation}
    where the first inequality follows from the triangle inequality, the second inequality follows from condition~\eqref{eq:cond} and the BRR error bound~\eqref{eq:brrerr}, the third inequality follows from the definition of $\delta$ in the component-wise context and the matrix norm inequality, the fourth inequality follows from condition~\eqref{eq:const1}, and the last equality follows from the definition of $\bar \alpha$.
\end{proof}

We can modify Proposition~\ref{prop:suffcond} to obtain an upper bound for the $\mathcal{V}$-norm of the error between the RB and HRBE solutions. To this end, we define $\lambda_{\min}$ and $\lambda_{\max}$ such that
\begin{equation}\label{eq:lambdas}
    \lambda_{\mn} = \inf_{v \in \mathcal{V}} \frac{\lr{v}_{\mathcal{V}}^2}{\lr{\mathbf{v}}_{2}^2},  \quad \lambda_{\mx} = \sup_{v \in \mathcal{V}} \frac{\lr{v}_{\mathcal{V}}^2}{\lr{\mathbf{v}}_{2}^2},
\end{equation}
and introduce the following corollaries.

\begin{corollary}[Absolute error bound]\label{coro:vboundabs}
    If all conditions of Proposition~\ref{prop:suffcond} hold, then
    \begin{equation}\label{eq:properrabsV}
    \begin{aligned}
        \lr{ {u}_{\rb}(\mu) - \widetilde{{u}}_{\rb}(\mu) }_{\mathcal{V}} &\leq (\bar{\alpha} +\bar{\varepsilon}) \sqrt{\lambda_\mx},
    \end{aligned}
    \end{equation}
    with the same $\bar{\alpha}$ and $\bar{\varepsilon}$ as in Proposition~\ref{prop:suffcond}.
\end{corollary}
\begin{proof}
    We first appeal to~\eqref{eq:lambdas} to obtain $\lr{{u}_{\rb}(\mu) - \widetilde{{u}}_{\rb}(\mu)}_{\mathcal{V}} \leq \sqrt{\lambda_\mx} \lr{\mathbf{u}_{\rb}(\mu) - \widetilde{\mathbf{u}}_{\rb}(\mu)}_2$. We then incorporate~\eqref{eq:properrabs} to obtain~\eqref{eq:properrabsV}.
\end{proof}

\begin{corollary}[Relative error bound]\label{coro:vboundrel}
    If all conditions of Proposition~\ref{prop:suffcond} hold and conditions~\eqref{eq:cond}, \eqref{eq:const1}, and \eqref{eq:lcond} are respectively replaced by
    \begin{align}
        \left\| \mathbf{u}_{\rb,c}(\mu) -  \bar{\mathbf{u}}_{\rb,c}(\mu) \right\|_2 &\leq  \bar{\varepsilon} \| {\mathbf{u}}_{\rb}(\mu) \|_2  &&\hspace{-4em} \forall c \in \sys, \nonumber\\
        \left\| \widetilde{\mathbf{R}}_{\rb,c}(\mathbf{P}^T_{\rb,c} \bar{\mathbf{u}}_{\rb}(\mu);\mu_c) \right\|_\infty &\leq \delta_{R_c} \| {\mathbf{u}}_{\rb}(\mu) \|_2  &&\hspace{-4em}\forall c \in \sys, \nonumber \\
        L \Big( 2 \| {\mathbf{u}}_{\rb}(\mu) \|_2 \frac{\sqrt{\sum_{c \in \sys} N_{M(c)} \delta_{R_c}^2 }} {\sigma - \sum_{c \in \sys} N_{M(c)} \delta_{J_c}} \Big) &\leq \frac{\sigma - \sum_{c \in \sys} N_{M(c)} \delta_{J_c}}{2 \sigma}, \nonumber
    \end{align}
    then
    \begin{equation}\label{eq:relerror}
        \begin{aligned}
        \frac{\left\| {u}_{\rb}(\mu) - \wtilde{{u}}_{\rb}(\mu) \right\|_{{\mathcal{V}}}}{\left\| {{u}}_{\rb}(\mu) \right\|_{{\mathcal{V}}}} &\leq (\bar{\alpha} + \bar{\varepsilon})\sqrt{\frac{\lambda_\max}{\lambda_\mn}},
        \end{aligned}
    \end{equation}
    with the same $\bar{\alpha}$ and $\bar{\varepsilon}$ as in Proposition~\ref{prop:suffcond}.
\end{corollary}
\begin{proof}
  We first observe
    \begin{equation}\label{eq:nmt1}\nonumber
        \begin{aligned}
            \lr{u_\rb(\mu) - \wtilde{u}_\rb(\mu)}_{\mathcal{V}} &\leq \lr{\mathbf{u}_\rb(\mu) - \wtilde{\mathbf{u}}_\rb(\mu)}_2 \sqrt{\lambda_\mx} \leq (\bar{\alpha} + \bar{\varepsilon}) \| {\mathbf{u}}_{\rb}(\mu) \|_2 \sqrt{\lambda_\mx},
        \end{aligned}
    \end{equation}
    where the first inequality follows from \eqref{eq:lambdas}, and the second equality follows from the application of Proposition~\ref{prop:suffcond} with $\bar{\varepsilon}$ and $\delta_{R_c}$, $c \in \sys$, replaced by $\bar{\varepsilon} \| {\mathbf{u}}_{\rb}(\mu) \|_2$ and $\delta_{R_c} \| {\mathbf{u}}_{\rb}(\mu) \|_2$, respectively. We finally appeal to \eqref{eq:lambdas} to obtain $\sqrt{\lambda_{\min}} \lr{{\mathbf{u}}_\rb(\mu)}_2 \leq \lr{{{u}}_\rb(\mu)}_{\mathcal{V}}$, which in turn yields \eqref{eq:relerror}.
\end{proof} 
\begin{remark}
  \label{rem:lambda_affine}
    The values of $\lambda_\mn$ and $\lambda_\mx$ are solely functions of the geometrical parameters in the system. The archetype components considered in this study, which will be introduced in Section~\ref{sec:res}, admit piecewise affine decompositions in their geometric parametrization. Consequently, computing $\lambda_\mn$ and $\lambda_\mx$ can be carried out efficiently in the online phase. Specifically, for each instantiated component $c \in \sys$, if $\{ {\phi}_{{c},i} \}_{i=1}^{\mathcal{N}_{{M(c)}}}$ denotes the geometric-parameter-dependent basis for ${\mathcal{V}}_{h,{c}}$ and $\mathbf{V}_c: \mathcal{D}_c \to \mathbb{R}^{\mathcal{N}_{M(c)} \times \mathcal{N}_{M(c)}}$ denotes its geometric-parameter-dependent inner-product matrix such that $(\mathbf{V}_c(\mu_c))_{i,j} = ({\phi}_{{c},j}, {\phi}_{{c},i})_{\mathcal{V}_c} \: \forall \mu_c \in \mathcal{D}_c$, $i,j = 1,\cdots, \mathcal{N}_{M(c)}$, the systems's geometric-parameter-dependent inner-product matrix $\mathbf{V}: \mathcal{D} \to \mathbb{R}^{\mathcal{N}_h \times \mathcal{N}_h}$ is given by $\mathbf{V}(\mu) = \sum_{c \in \sys} \mathbf{P}_{c} \mathbf{V}_c(\mu_c) \mathbf{P}_{c}^T \: \forall \mu \in \mathcal{D}$. Here, $\mathbf{P}_c: \mathbb{R}^{\mathcal{N}_{M(c)}} \to \mathbb{R}^{\mathcal{N}_h} \: \forall c \in \sys$ are linear extension operators that map the components truth to system's DoF. Since (we have assumed) $\mathbf{V}_c(\cdot) \: \forall c \in \sys$ admit piecewise affine decompositions, forming $\mathbf{V}_c(\cdot)$ and hence $\mathbf{V}(\cdot)$ during the online phase does not rely on the components' truth FE discretizations and quadrature rules. Additionally, the computation of extreme eigenvalues of $\mathbf{V}(\cdot)$ (i.e., $\lambda_\mn$ and $\lambda_\mx$) can be performed efficiently using iterative methods such as the Lanczos algorithm~\cite{trefethen2022numerical}.
\end{remark}

\subsection{Component-wise hyperreduction: formulation}
\label{subsec:eqptraining}

Using Proposition~\ref{prop:suffcond} and Corollaries~\ref{coro:vboundabs} and \ref{coro:vboundrel}, we now develop a component-wise hyperreduction training routine for the archetype components in the library. In Section~\ref{subsec:train_set_gen}, we introduced for each archetype component $\wc \in \lib$, a training parameter set $\Xi_\wc^\train$ and its corresponding state training set $U_{h,\wc}^{\train}$ (Algorithm~\ref{alg:emptr}). In Section~\ref{subsec:rbptraining}, we also described a procedure to construct an RB for its bubble space $\widehat{\mathcal{V}}_{\rb,\wc}^\bb$ using its associated bubble training set $U_{h,\wc}^{\train, \bb}$. Since hyperreduction is carried out with respect to the RB solutions, for each archetype component $\wc$, we define a state training set $U_{\rb,\wc}^{\train} \equiv \{ u_{\rb,\wc,n}^{\train}\}_{n=1}^{N_\sample}$ (or its algebraic equivalent $\mathbf{U}_{\rb,\wc}^{\train} \equiv \{ \mathbf{u}_{\rb,\wc,n}^{\train}\}_{n=1}^{N_\sample}$), where the RB snapshots $u_{\rb,\wc,n}^{\train}, \: n \in \{1, \cdots, N_\sample\} $, are generated using Algorithm~\ref{alg:emptreqp}. 

\begin{algorithm}[!tb]
    \caption{Generating bubble RB snapshots for hyperreduction of archetype components.}\label{alg:emptreqp}
    \SetAlgoLined
    \KwInput{The previously generated set of truth snapshot solutions $U_{h,\wc}^{\train}$ $\forall \wc \in \lib$}
    \KwOutput{The set of RB snapshot solutions $U_{\rb,\wc}^{\train}$ $\forall \wc \in \lib$}
        \For{$\wc \in \lib$}{
            $U_{\rb,\wc}^{\train} = \emptyset$\;
            \For{$n=1, \cdots, N_{\sample}$}{
                Decompose $u_{h,\wc,n}^{\train}$ in the previously constructed training set $U_{h,\wc}^{\train}$ in Algorithm~\ref{alg:emptr} into bubble $u_{h,\wc,n}^{\train,\bb}$ and port solutions $\{ u_{h,\wc,n}^{\train,\wp} \}_{\wp \in \pset_\wc}$, as in~\eqref{eq:compsol}\;
                Compute the restriction of port solutions on the ports (i.e., $\{ u_{h,\wc,n}^{\train,\wp}\Big|_{\gamma_{\wc,\wp}} \}_{\wp \in \pset_\wc}$)\;  		
                Compute $u_{\rb,\wc,n}^{\train,\bb}$ by solving~\eqref{eq:rbtruth} for a system composed of only component $\wc$ with $\{ u_{h,\wc,n}^{\train,\wp}\Big|_{\gamma_{\wc,\wp}} \}_{\wp \in \pset_\wc}$ as the Dirichlet boundary conditions imposed on its $n^\gamma_\wc$ ports\;
                Compute $u_{\rb,\wc,n}^{\train} = u_{\rb,\wc,n}^{\train,\bb} + \sum_{\wp \in \pset_\wc}u_{h,\wc,n}^{\train,\wp}$\;
                $U_{\rb,\wc}^{\train} \leftarrow U_{\rb,\wc}^{\train}  \cup u_{\rb,\wc,n}^{\train}$\;
            }
        }
\end{algorithm}

Additionally, for each archetype component $\wc \in \lib$, we introduce the algebraic RB residual and Jacobian ${\mathbf{R}}_{\rb,\wc}: \mathbb{R}^{N_\wc} \times \what{\mathcal{D}}_\wc \rightarrow \mathbb{R}^{N_{\wc}}$ and ${\mathbf{J}}_{\rb,\wc}: \mathbb{R}^{N_{\wc}} \times \what{\mathcal{D}}_\wc \rightarrow \mathbb{R}^{N_{\wc} \times N_{\wc}}$ formulated in Appendix~\ref{app:resandjacrb}. We further denote the \emph{barred} versions of the introduced RB algebraic terms. These barred versions are formulated the same as their respective RB counterparts, albeit with the truth quadrature weights $\{ \what{{\rho}}_{\wc,q} \}_{q=1}^{Q_{\wc}} \: \: \forall \wc \in \lib$ in \eqref{eq:resbrb} replaced by $\bar{\rho}_\wc \equiv \{ {\bar{\rho}}_{\wc,q} \}_{q=1}^{Q_{\wc}} \:\: \forall \wc \in \lib$, which are the design variables (unknowns) for the hyperreduction problem. Then, we pose the component-wise hyperreduction problem in the offline phase for $\wc \in \lib$ as follows: given a parameter training set $\Xi^\train_\wc$, state training set $U_{\rb,\wc}^{\train}$ (or its algebraic equivalent $\mathbf{U}_{\rb,\wc}^\train$), domain volume $|\Omega_\wc|$, and hyperparameter ${\delta_{\wc}}$, find ${\bar{\rho}}_\wc^* \in\mathbb{R}^{Q_\wc}$ such that
\begin{alignat}{3}\label{eq:finaleqp}
  \bar{\rho}_{\wc}^* = \argmin_{  \{\bar{\rho}_{{\wc},q}\}_{q=1}^{Q_{{\wc}}} } \quad \|  \bar{\rho}_{{\wc},q} \|_0
\end{alignat}
subject to
\begingroup
\allowdisplaybreaks
\begin{alignat}{4}
\bar{\rho}_{{\wc},q} &\geq 0,  \quadd &&q=1,\cdots,Q_{\wc}, \\
\Big| |\Omega_{\wc}| - \sum_{q=1}^{Q_{\wc}} \bar{\rho}_{{\wc},q} \Big| &\leq \delta_{\wc}, \label{eq:compconst01}\\
\left\| \bar{\mathbf{R}}_{\rb,{\wc}}({\mathbf{u}}^\train_{\rb,{\wc}}(\mu);\mu, \bar{\rho}_{\wc}) \right\|_\infty &\leq \delta_{\wc} \quadd  &&\forall \mu \in \Xi^\train_{\wc},\label{eq:compconst11}\\
\left\| \mathbf{J}_{\rb,{\wc}}(\mathbf{u}^\train_{\rb,{\wc}}(\mu); \mu)-  \bar{\mathbf{J}}_{\rb,{\wc}}({\mathbf{u}}^\train_{\rb,{\wc}}(\mu);\mu, \bar{\rho}_{\wc}) \right\|_\mx &\leq  \delta_{\wc} \quadd  &&\forall \mu \in \Xi^\train_{\wc}. \label{eq:compconst21}
\end{alignat}
\endgroup
The $\ell^0$-minimization problem seeks the sparsest quadrature rule that satisfies the constraints. In practice, we approximate the $\ell^0$-minimization problem as an $\ell^1$-minimization problem (with the objective function $\sum_{q=1}^{Q_{{\wc}}} \bar{\rho}_{{\wc},q}$) and solve the problem using a simplex method following \cite{yano2019lp}. The enforcement of the constant function constraint~\eqref{eq:compconst01} enhances the robustness of the hyperreduction training and is a reasonable condition for any quadrature scheme \cite{yano2019lp}. The RQ rule for each component is determined by $( \th{x}_{{\wc},q}^r , \th{\rho}_{{\wc},q}^r )_{q=1}^{\wtilde{Q}_{\wc}^r} = ( ( \widehat{x}_{{\wc},q}, \bar{\rho}^*_{{\wc},q} ) | \: \bar{\rho}^*_{{\wc},q} > 0 )_{q=1}^{Q_{\wc}}$.  In this work, we set $\delta_{R_\wc} = \delta_{J_\wc} = \delta_{\wc}$. To construct the RQ rule $( \th{x}_{{\wc},q}^f , \th{\rho}_{{\wc},q}^f )_{q=1}^{\wtilde{Q}_{\wc}^f}$ for the output~\eqref{eq:outputnlrbe}, we follow the procedure in \cite{yano2019lp} and replace the constraints~\eqref{eq:compconst01}--\eqref{eq:compconst21} with analogous constraints for the output functional $F(\cdot; \mu)$ and solve the EQP optimization problem.

\section{Offline and online computational procedure}
\label{sec:comppro}
In this section, we develop the offline--online computational procedure for the HRBE method. A key challenge to offline--online computational decomposition that provides quantitative control of the hyperreduction error at the system level is this: the hyperreduction training is performed for each archetype component independently in the offline phase; therefore, unlike in the monodomain setting for which the EQP is originally designed (e.g.,~\cite{yano2019lp}), the minimum singular value of the Jacobian of the ultimate systems created by assembling the trained archetype components, which is required in~\eqref{eq:properrabs}, \eqref{eq:properrabsV}, and \eqref{eq:relerror} to control the error, is not available at the training time. To address this challenge, we propose an approach where the hyperreduction training for any $\wc \in \lib$ is conducted in the offline phase with various ${\delta}_{\wc}$ values. Subsequently, in the online phase, the appropriate RQ rule for each component is adaptively chosen and applied to solve the HRBE problem through an iterative bootstrap process. We now present the offline--online computational procedure.

\subsection{Offline phase}
\label{subsec:offcost}
In the offline stage, we prepare the RB $\{ \widehat{\xi}_{\wc,i}^\bb\}_{i=1}^{N_{\wc}^\bb}$ of the bubble spaces $\widehat{\mathcal{V}}_{\rb,\wc}^\bb$ and the RQ rules $(\th{x}_{\wc,q}^r, \th{\rho}_{\wc,q}^r )_{q=1}^{\wtilde{Q}^r_\wc}$ and $(\th{x}_{\wc,q}^f, \th{\rho}_{\wc,q}^f )_{q=1}^{\wtilde{Q}^f_\wc}$ for each of $N_\arch$ archetype components $\wc \in \lib$. To construct the RB, we first use Algorithm~\ref{alg:emptr} to generate the training set $U^\train_{h,\wc}$ for each archetype component $\wc$. For the $n$-th sample subsystem, $n \in \{1,\cdots,N_\sample\}$, of the archetype component $\wc$, the computation of the solution $u_{h,\wc,n}^{\train}$ in Line~\ref{ln:solve} of Algorithm~\ref{alg:emptr} requires $\mathcal{O}({Q}_{h,\wc,n}^\sub)$ operations for the assembly of the residual and Jacobian and $\mathcal{O}((\mathcal{N}_{h,\wc,n}^\sub)^l)$ operations for the solution of the linear system per Newton iteration, where $Q_{h,\wc,n}^\sub$ is the subsystem's number of truth quadrature points, $\mathcal{N}_{h,\wc,n}^\sub$ is the subsystem's number of truth DoF, and the coefficient $1 \leq l \leq 2$ depends on the solver and the domain dimension. Typical problems that we consider require 5 to 15 Newton iterations for convergence. The subsequent computational cost of the POD is negligible compared to the cost to generate the training set.

We now analyze the cost of hyperreduction for each archetype component $\wc \in \lib$. Using Algorithm~\ref{alg:emptreqp} to generate the RB snapshots requires the solution of a nonlinear system of equations of size $N_\wc^\bb$ for each training sample. This incurs, for each snapshot, a cost of $\mathcal{O}(N_\wc^2 Q_{\wc})$ operations for computing the residual and Jacobian and a cost of $\mathcal{O}((N_\wc^\bb)^3)$ operations for solving the linear system in each Newton iteration. Additionally, computing the outputs needed in output hyperreduction of each component requires $\mathcal{O}(Q_\wc)$ operations for each training sample. Moreover, for each archetype component, a simplex method is used to approximately solve the hyperreduction problem~\eqref{eq:finaleqp}--\eqref{eq:compconst21} for different $\delta_\wc$ values. Each problem has $Q_\wc$ unknowns, $Q_\wc$ positivity constraints, 1 constant function constraint, $N_\sample N_{\wc}$ residual constraints, and $N_\sample N_{\wc}^2$ Jacobian constraints. In practice, the absolute value constant function, residual, and Jacobian constraints are converted into $2 (1 + N_\sample N_{\wc} + N_\sample N_{\wc}^2)$ inequality constraints. Additionally, the output hyperreduction involves the solution of an optimization problem with $Q_\wc$ unknowns and $2 (1 + N_\sample)$ inequality constraints.

\subsection{Online phase: adaptive RQ selection}
\label{subsec:nlrbecompr}

We now describe a procedure to find the RQ rule of each component in the system in the online phase such that, for any given topological configuration and $\mu \in \mathcal{D}$, the HRBE solution $\wtilde{{u}}_\rb(\mu)$ achieves the target $\mathcal{V}$-norm error with respect to the RB solution ${{u}}_\rb(\mu)$. Our formulation builds on Corollaries~\ref{coro:vboundabs} and \ref{coro:vboundrel}. We note that since the error is measured with respect to ${u}_\rb(\mu)$, it is implicitly assumed that $\bar{\mathbf{u}}_\rb(\mu) = \mathbf{u}_\rb(\mu)$ in Proposition~\ref{prop:suffcond}, and $\bar{\varepsilon} = 0$ in~\eqref{eq:properrabs}, \eqref{eq:properrabsV}, and \eqref{eq:relerror}.

We first discuss an online-efficient procedure to compute $\sigma_{\mn}(\mathbf{J}_{\rb}\left(\mathbf{u}_\rb(\mu);\mu) \right)$ (or more precisely approximate it), required for computing $\bar{\alpha}$ in~\eqref{eq:properrabs}, \eqref{eq:properrabsV}, and \eqref{eq:relerror}. A direct computation of $\sigma_{\mn}(\mathbf{J}_{\rb}\left(\mathbf{u}_\rb(\mu);\mu) \right)$ poses two computational challenges. Firstly, the computation requires the RB solution $\mathbf{u}_\rb(\mu)$, which defeats the purpose of hyperreduction; we wish to use only its HRBE counterpart $\wtilde{\mathbf{u}}_\rb(\mu)$. Secondly, it involves forming the RB Jacobian $\mathbf{J}_{\rb}(\cdot;\mu)$, which depends on the components' truth quadrature rules and prevents efficient online computation.

To address these challenges, we appeal to the BRR theorem. We set $G(\cdot) \equiv \mathbf{R}_\rb(\cdot; \mu)$ and $\mathbf{v} \equiv \wtilde{\mathbf{u}}_\rb(\mu)$ in Lemma~\ref{lemma:brr}, assume $\mathbf{J}_\rb(\wtilde{\mathbf{u}}_\rb(\mu); \mu)$ is nonsingular and the conditions of the theorem hold, and apply~\eqref{eq:brrjac} to obtain
\begin{equation}\label{eq:sigmaest}
\begin{aligned}
\lr{\mathbf{J}_\rb^{-1}({\mathbf{u}}_\rb(\mu); \mu)}_2 = \frac{1}{\sigma_{\mn}(\mathbf{J}_{\rb}\left({\mathbf{u}}_{\rb}(\mu);\mu) \right) } \leq 2 \lr{\mathbf{J}_\rb^{-1}(\wtilde{\mathbf{u}}_\rb(\mu); \mu)}_2 = \frac{2}{\sigma_{\mn}(\mathbf{J}_{\rb}\left(\wtilde{\mathbf{u}}_{\rb}(\mu);\mu) \right) }.
\end{aligned}
\end{equation}
Therefore, $\sigma_{\mn}(\mathbf{J}_{\rb}\left(\wtilde{\mathbf{u}}_{\rb}(\mu);\mu) \right)/2$ is a lower bound for $\sigma_{\mn}(\mathbf{J}_{\rb}\left({\mathbf{u}}_{\rb}(\mu);\mu) \right)$. In order to approximate $\sigma_{\mn}(\mathbf{J}_{\rb}\left(\wtilde{\mathbf{u}}_{\rb}(\mu);\mu) \right)$ we appeal to the following lemma.

\begin{lemma}\label{lemma:singularpert}
For any three matrices $\mathbf{A} \in \mathbb{R}^{N \times N}$, $\mathbf{B} \in \mathbb{R}^{N \times N}$, and $\mathbf{C} \in \mathbb{R}^{N \times N}$ such that $\mathbf{A} = \mathbf{B} + \mathbf{C}$
\begin{equation}\label{eq:singularpert}
| \sigma_\mn(\mathbf{A}) - \sigma_\mn(\mathbf{B}) | \leq \sigma_\mx(\mathbf{C}),
\end{equation}
where $\sigma_\mn (\cdot)$ and $\sigma_\mx(\cdot)$, respectively, correspond to the minimum and maximum singular values of their argument.
\end{lemma}
\begin{proof}
  We first observe that, for $\mathbf{A} = \mathbf{B} + \mathbf{C}$,
\begin{equation}\nonumber
\begin{aligned}
\sigma_\mn (\mathbf{A}) = \min_{\mathbf{v} \in \mathbb{R}^N}  \frac{ \lr{(\mathbf{B} + \mathbf{C}) \mathbf{v}}_2 }{\lr{\mathbf{v} }_2}  \geq \min_{\mathbf{v} \in \mathbb{R}^N}  \frac{ \lr{\mathbf{B} \mathbf{v}}_2 - \lr{\mathbf{C} \mathbf{v}}_2 }{\lr{\mathbf{v} }_2} &\geq \min_{\mathbf{v} \in \mathbb{R}^N}  \frac{ \lr{\mathbf{B} \mathbf{v}}_2}{\lr{\mathbf{v} }_2}  - \max_{\mathbf{v} \in \mathbb{R}^N}  \frac{ \lr{\mathbf{C} \mathbf{v}}_2 }{\lr{\mathbf{v} }_2} \\
&= \sigma_\mn(\mathbf{B}) - \sigma_\mx(\mathbf{C})
\end{aligned}
\end{equation}
and hence  $\sigma_\mn(\mathbf{B}) - \sigma_\mn (\mathbf{A}) \leq \sigma_\mx(\mathbf{C})$, 
where the first and last equality follow from the definition of the extreme singular values, and the first inequality follows from the triangle inequality. We apply an analogous sequence of inequalities to $\mathbf{B} = \mathbf{A} - \mathbf{C}$ to obtain $\sigma_\mn(\mathbf{A}) - \sigma_\mn (\mathbf{B}) \leq \sigma_\mx(\mathbf{C})$. The combination of the two inequalities yields the desired result.
\end{proof}
The application of the lemma to $\mathbf{J}_\rb(\wtilde{\mathbf{u}}_{\rb}(\mu); \mu) = \widetilde{\mathbf{J}}_\rb(\wtilde{\mathbf{u}}_{\rb}(\mu); \mu) + \pr{\mathbf{J}_\rb(\wtilde{\mathbf{u}}_{\rb}(\mu); \mu) -  \widetilde{\mathbf{J}}_\rb(\wtilde{\mathbf{u}}_{\rb}(\mu); \mu)}$ yields
\begin{equation}\nonumber
| \sigma_\mn(\mathbf{J}_\rb(\wtilde{\mathbf{u}}_{\rb}(\mu); \mu)) - \sigma_\mn( \widetilde{\mathbf{J}}_\rb(\wtilde{\mathbf{u}}_{\rb}(\mu); \mu)) | \leq \sigma_\mx(\mathbf{J}_\rb(\wtilde{\mathbf{u}}_{\rb}(\mu); \mu)- \widetilde{\mathbf{J}}_\rb(\wtilde{\mathbf{u}}_{\rb}(\mu); \mu)).
\end{equation}
In other words, as the disparity between the RB and hyperreduced RB Jacobians decreases, so does the discrepancy between $\sigma_\mn(\mathbf{J}_\rb(\wtilde{\mathbf{u}}_{\rb}(\mu); \mu))$ and $\sigma_\mn( \widetilde{\mathbf{J}}_\rb(\wtilde{\mathbf{u}}_{\rb}(\mu); \mu))$. Given that the hyperreduction training for each archetype component is intended to reduce this very gap, we propose to use $\sigma_\mn( \widetilde{\mathbf{J}}_\rb(\wtilde{\mathbf{u}}_{\rb}(\mu); \mu))$ in place of $\sigma_\mn(\mathbf{J}_\rb(\wtilde{\mathbf{u}}_{\rb}(\mu); \mu))$. Finally, we combine this approximation with the lower-bound estimate~\eqref{eq:sigmaest} to conservatively approximate $\sigma_{\mn}(\mathbf{J}_{\rb}\left(\mathbf{u}_\rb(\mu);\mu) \right)$ by $\sigma_{\mn}(\wtilde{\mathbf{J}}_{\rb}\left(\wtilde{\mathbf{u}}_{\rb}(\mu);\mu) \right) / 2$.

We finally propose the adaptive procedure, Algorithm~\ref{alg:sing}, to find the components' RQ rules and the HRBE solution in the online phase. For a given $\mu \in \mathcal{D}$ and a desired system-level $\mathcal{V}$-norm error $\epsilon$ between the RB and HRBE solutions, we first compute $\bar{\alpha} = \epsilon / \sqrt{\lambda_{\mx}}$ for absolute error control (or $\bar{\alpha} = \epsilon / \sqrt{\lambda_{\mx} / \lambda_\mn}$ for relative error control). We then use the RQ rules associated with the initial $\delta_{c}$ values $\forall c \in \sys$ to compute the HRBE solution $\wtilde{\mathbf{u}}_{\rb}(\mu)$ and $\sigma \equiv \sigma_\mn( \widetilde{\mathbf{J}}_\rb(\wtilde{\mathbf{u}}_{\rb}(\mu); \mu))/2$. Then, for all $c \in \sys$, we set ${\delta}_{R_c} = {\delta}_{J_c} =  {\delta}_{c}$ for absolute error control (or ${\delta}_{R_c}= {\delta}_{J_c} = {\delta}_{c} / \lr{\wtilde{\mathbf{u}}_{\rb}(\mu)}_2$ for relative error control). If $\sum_{c \in \sys} N_{M(c)} {\delta}_{J_c} \geq \sigma$ or $2{{\sum_{c \in \sys} \sqrt{N_{M(c)}} {\delta}_{R_c}}} / (\sigma - \sum_{c \in \sys} N_{M(c)} {\delta}_{J_c}) > \bar{\alpha}$, the hyperreduction tolerances ${\delta}_{c} = {\delta}_{R_c} = {\delta}_{J_c}$ of each component $c \in \sys$ is adjusted such that these conditions hold. We then use the  RQ rules associated with the new hyperreduction tolerances to compute the new HRBE solution $\wtilde{\mathbf{u}}_{\rb}(\mu)$ and $\sigma = \sigma_{\mn}(\wtilde{\mathbf{J}}_{\rb}\left(\wtilde{\mathbf{u}}_{\rb}(\mu);\mu) \right)/2$. This process is repeated until convergence; for the problems considered in Section~\ref{sec:res}, the procedure converges in two iterations.

\begin{algorithm}[!tb]
\caption{Adaptive selection of RQ rules and solving the HRBE problem in the online phase.}\label{alg:sing}
\SetAlgoLined
\KwInput{System-level $\mu \in \mathcal{D}$ and desired $\mathcal{V}$-norm error $\epsilon \in \mathbb{R}_{>0}$ between RB and HRBE solutions}
\KwOutput{The HRBE solution and components' RQ rules}
    Compute $\lambda_\mx$ (and $\lambda_\mn$ if $\epsilon$ is the relative error) for the system\;
    Set $\bar{\alpha} = \epsilon / \sqrt{\lambda_\mx}$ (or $\bar{\alpha} = \epsilon / \sqrt{\lambda_\mx / \lambda_\mn}$ if $\epsilon$ is the relative error)\;
	Select the initial hyperreduction tolerances $\delta_{c}$ $\forall c \in \sys$\;
	\While{true}{\label{ln:sing3}
        Set the RQ rules associated with the current $\delta_{c}$ values $\forall c \in \sys$\;\label{ln:sing1}
		Solve the HRBE problem to find $\wtilde{\mathbf{u}}_\rb(\mu)$\;
		Find $\sigma \equiv \sigma_\mn( \widetilde{\mathbf{J}}_\rb(\wtilde{\mathbf{u}}_{\rb}(\mu); \mu))/2$\;
                Set ${\delta}_{R_c} = {\delta}_{J_c} = {\delta}_{c}$ for all $c \in \sys$ (or ${\delta}_{R_c} = {\delta}_{J_c} = {\delta}_{c} / \lr{\wtilde{\mathbf{u}}_{\rb}(\mu)}_2$ if $\epsilon$ is the relative error)\;
		\eIf{$\sum_{c \in \sys} N_{M(c)} {\delta}_{J_c} \geq \sigma$ or $2{{\sum_{c \in \sys} \sqrt{N_{M(c)}} {\delta}_{R_c}}} / (\sigma - \sum_{c \in \sys} N_{M(c)} {\delta}_{J_c}) > \bar{\alpha}$}{\label{ln:convcond}
            Update ${\delta}_{c}$ and subsequently ${\delta}_{R_c}$ and $\delta_{J_c}$ $\forall c \in \sys$ such that both conditions hold\;
			Go to Step~\ref{ln:sing1}\;
		}{
			break\;
		}
	}\label{ln:sing4}
\end{algorithm}

\begin{remark}
  \label{rem:no_adaptive_rb}
  In this work, we do \emph{not} consider adaptive selection of the RB for each component to control the truth vs RB error. Instead, we consider an adaptive selection of hyperreduction tolerance $\delta_c$, and hence the RQ rules, for each component to achieve the desired system-level hyperreduction error. Thus, the RB is fixed independent of the target hyperreduction tolerance for each archetype component. We focus on developing online-adaptive hyperreduction for component-based systems, and defer the development of online-adaptive RB selection for component-based systems to future work.
\end{remark}

\subsection{Online phase: computational cost and memory footprint}

We now remark on the computational cost of solving the HRBE problem using Algorithm~\ref{alg:sing} as opposed to that for solving the truth problem. For the truth problem, each iteration of Newton's method necessitates $\mathcal{O}(Q_h \equiv \sum_{c \in \sys} Q_{M(c)})$ operations for evaluating the truth residual and Jacobian. Additionally, the solution of the linear system~\eqref{eq:newtonupdate} in every Newton iteration requires $\mathcal{O}(\mathcal{N}_h^n)$ operations, where $1 \leq n \leq 2$ depends on the domain dimension and the solver employed. Moreover, computing the truth output involves $\mathcal{O}(Q_h)$ operations.

On the other hand, each cycle of the loop in Algorithm~\ref{alg:sing} (Lines~\ref{ln:sing3}--\ref{ln:sing4}) involves solving the HRBE problem and computing $\sigma_\mn( \widetilde{\mathbf{J}}_\rb(\wtilde{\mathbf{u}}_{\rb}(\mu); \mu))$. In each Newton iteration, evaluating the HRBE residual $\widetilde{\mathbf{R}}_{\rb}(\cdot;\mu)$ and Jacobian $\widetilde{\mathbf{J}}_{\rb}(\cdot;\mu)$ requires $\mathcal{O}(\sum_{c \in \sys} N_{\rb}^2 \wtilde{Q}_{M(c)}^r) \ll \mathcal{O}(Q_h)$ operations, where $\wtilde{Q}_{\wc}^r$ is the number of RQ points of the archetype component $\wc \in \lib$ in a given cycle. In addition, finding the Newton update requires solving a linear system---which is component-block-wise sparse---in $\mathcal{O}(N_{\rb}^n) \ll \mathcal{O}(\mathcal{N}_h^n)$ operations. Also, computing the minimum singular value involves $\mathcal{O}(N_{\rb}^n)$ operations. Finally, once $\wtilde{\mathbf{u}}_{\rb}(\mu)$ is found, computing the approximate output $\widetilde{F}_{\rb}(\cdot;\mu)$ requires $\mathcal{O}(\sum_{c \in \sys} \wtilde{Q}_{M(c)}^f) \ll \mathcal{O}(Q_h)$ operations.

We now compare the memory footprint of the truth and HRBE formulations. The storage requirement for the truth problem, dominated by the truth Jacobian storage, is $\mathcal{O}\left( \mathcal{N}_h^n \right)$; $n = 1$ if an iterative solver is used at each iteration of the Newton method, otherwise $n = 4/3$ for $d \leq 3$ to store factorization. For the HRBE problem, the entire library must be loaded in the computer memory. To compute the residual $\widetilde{\mathbf{R}}_{\rb}(\cdot;\mu)$, Jacobian $\widetilde{\mathbf{J}}_{\rb}(\cdot;\mu)$, and output functional $\widetilde{F}(\cdot;\mu)$, we precompute and store the following quantities for each archetype component $\wc \in \lib$: (i) the RQ rules $( \th{x}_{\wc,q}^r, \th{\rho}_{\wc,q}^r )_{q=1}^{\wtilde{Q}_\wc^r}$ and $( \th{x}^f_{\wc,q} , \th{\rho}^f_{\wc,q} )_{q=1}^{\wtilde{Q}_\wc^f}$ for different $\delta_\wc$ values, (ii) the values of the bubble space basis $\{ \widehat{\xi}_{\wc,i}^\bb \}_{i=1}^{N_\wc^\bb}$ and the port basis $\{ \widehat{\psi}_{\wc,i}^\wp \}_{i=1}^{\mathcal{N}_{\wc}^\wp}$ for all ports $\wp \in \pset_\wc$ at the RQ points $\{\th{x}_{\wc,q}^r \}_{q=1}^{\wtilde{Q}_\wc^r}$ and $\{\th{x}^f_{\wc,q} \}_{q=1}^{\wtilde{Q}_\wc^f}$ associated with different $\delta_\wc$ values, and (iii)  the gradient values of the bubble space basis $\{ \nabla \widehat{\xi}_{\wc,i}^\bb \}_{i=1}^{N_\wc^b}$ and the port basis $\{ \nabla \widehat{\psi}_{\wc,i}^\wp \}_{i=1}^{\mathcal{N}_{\wc}^\wp}$ for all ports $\wp \in \pset_\wc$ at the RQ points $\{\th{x}_{\wc,q}^r \}_{q=1}^{\wtilde{Q}_\wc^r}$ and $\{\th{x}^f_{\wc,q} \}_{q=1}^{\wtilde{Q}_\wc^f}$ associated with different $\delta_\wc$ values. Therefore, the total online storage for $\wc \in \lib$ is
\begin{equation}\nonumber
(d + 1) \sum_{\wc \in \lib} N_{\delta_\wc} \left( \wtilde{Q}_\wc^r + \wtilde{Q}_\wc^f \right)  \left( 1 + N_\wc^\bb +  \sum_{\wp \in \pset_\wc} \mathcal{N}_{\wc}^\wp \right),
\end{equation}
where $N_{\delta_\wc}$ is the number of hyperreduction tolerances of $\wc \in \lib$ for which the state and output hyperreduction trainings are performed. Therefore, the online storage is independent of $\mathcal{N}_{\wc}^{\bb}$ and $Q_\wc$ $\forall \wc \in \lib$. Furthermore, owing to $N_\wc^\bb \ll \mathcal{N}_{\wc}^\bb$, $\wtilde{Q}_\wc^r \ll Q_\wc$, and $\wtilde{Q}_\wc^f \ll Q_\wc \: \: \forall \wc \in \lib$, the online storage requirement for the HRBE problem is significantly smaller than that of the truth problem, particularly when the truth DoF is large and there is a significant reuse of archetype components: i.e., $N_\arch$ is small relative to the size of the system, which is the case for which the HRBE method is designed. It is important to note that the storage requirement scales with $N_\arch$ rather than ${N_\comp}$. Therefore, employing the HRBE method ensures that the storage and computational cost of the online phase are independent of $\mathcal{N}_h$ and $Q_h$, as desired.

\section{Example: nonlinear thermal fin systems}
\label{sec:res}
\subsection{Problem description}
\label{subsec:description}
We now apply the HRBE method to two-dimensional nonlinear thermal fin systems. Systems are made of an aluminum alloy \cite{nist_aluminum_3003} with a temperature-dependent thermal conductivity $k: [1, 300] \: \mathrm{K} \to [4.341, 177.868] \: \mathrm{W/K}$ that satisfies
\begin{equation}\label{eq:conductivity}
\log (k(x)) = \sum_{i=0}^7 k_i \: (\log (x))^i \quad \forall x \in [1, 300] \: \mathrm{K},
\end{equation}
where $k_i$, $i=1,\cdots,7$, are given in Table~\ref{table:coeffs}. The parameterized continuous residual form for the ultimate systems is
\begin{equation}\nonumber
R(w,v;\mu) =  \int_{\Omega(\mu)} \left(k(w) \nabla w \right) \cdot \nabla v \: dx - \int_{\Omega(\mu)} f(\mu) \: v \: dx \quad \forall w,v \in \mathcal{V},
\end{equation} 
where $\mathcal{V} \equiv \left\{ v \in H^1(\Omega(\mu)) \Big| \ v_{\Gamma_D} = 0 \right\}$, and $f: \mathcal{D} \to L^2(\Omega(\mu))$ is the volumetric source term, which is assumed to be constant within each component. The residual form does not admit an affine decomposition, as the first integral depends nonlinearly on the field variable.

\begin{remark}
  \label{rem:hrbe_generality}
  The HRBE method, like the (monodomain) RB-EQP method \cite{yano2019lp}, can also treat other types of nonlinear solution and parameter dependencies, including those arising from nonlinear geometric transformations from reference to physical domains of the instantiated components.
\end{remark}

\begin{table}
\caption{\label{table:coeffs}Coefficients of the aluminum's thermal conductivity equation~\eqref{eq:conductivity}.}
\begin{center}
\begin{tabular}{lcccccccc}
\toprule[0.1em]
{Coefficient} & $k_0$ & $k_1$ & $k_2$ & $k_3$ & $k_4$ & $k_5$ & $k_6$ & $k_7$\\[0.2em]
{Value (W/K)} & 0.637 & -1.144 & 7.462 & -12.691 & 11.917 & -6.187 & 1.639 & -0.173\\[0.2em]
\bottomrule
\end{tabular}
\end{center}
\end{table}

\subsection{Archetype component library}
\label{subsec:library}
\begin{figure}
\centering
\includegraphics[width=1\textwidth]{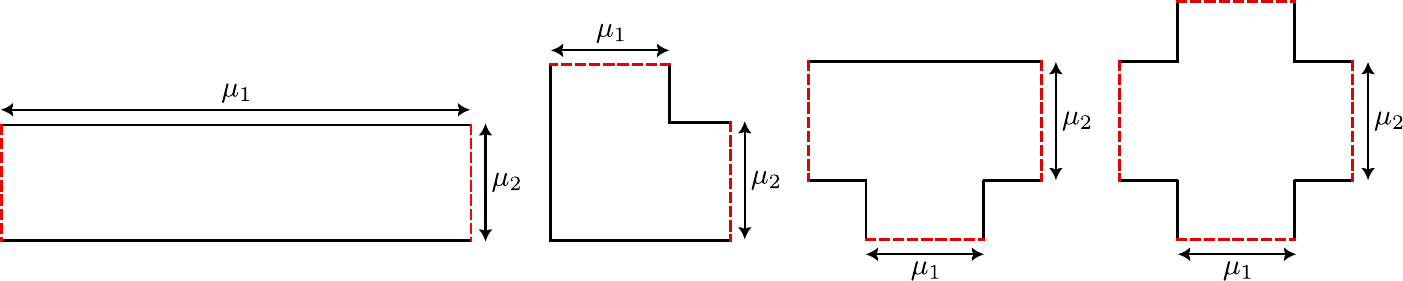}
        \caption{Archetype components in their reference domains. From left to right: rod, bracket, tee and cross. Local ports are shown by red dashed lines.}\label{fig:archcomps}
\end{figure}

Our archetype component library comprises four archetype components as shown in Figure~\ref{fig:archcomps}. Each archetype component is characterized by two geometric parameters $\mu_1$ and $\mu_2$, and one physical parameter $\mu_3 \in [0, 10]$ $\mathrm{W/cm^2}$ associated with volumetric heat source. For all components, $\mu_1 \in [0.5, 1]$ cm and $\mu_2 \in [0.5, 1]$ cm, with the exception of the rod component where $\mu_1 \in [3, 6]$ cm. The values of geometric parameters $\mu_1$ and $\mu_2$ in the reference domain of all archetype components are 1 cm, except for $\mu_1$ of the rod component, which is 4 cm. All components admit piecewise affine geometric transformations from their reference to physical spatial domains. Therefore, $\lambda_\mn$ and $\lambda_\mx$ of the systems instantiated from these components, required in Algorithm~\ref{alg:sing}, can be computed efficiently during the online phase; cf.\ Remark~\ref{rem:lambda_affine}. As shown in Figure~\ref{fig:archcomps}, the rod and bracket components have two local ports, the tee component has three local ports, and the cross component has four local ports. All ports are mapped from the same 17-DoF reference port discretized by eight quadratic line elements. Furthermore, all components are discretized using quadratic triangular elements leading to $\mathcal{N}^{\bb}_{\mathrm{rod}} = 691$, $\mathcal{N}^{\bb}_{\mathrm{bracket}} = 703$, $\mathcal{N}^{\bb}_{\mathrm{tee}} = 1026$, and $\mathcal{N}^{\bb}_{\mathrm{cross}} = 1165$.

The offline training proceeds in three sequential steps. First, for each archetype component we generate a set of empirical training data using Algorithm~\ref{alg:emptr}. Specifically, for each target component, we create $N_{\sample} = 100$ sample subsystems by connecting it with a probability of $\beta = 0.8$ to other components in the library through each of its ports. We then assign uniformly random parameter values to the components in the subsystems and set uniformly random constant Dirichlet boundary conditions to their nonshared global ports, ranging from $1$~K to $250$~K.

\begin{figure}
\centering
    \begin{subfigure}{0.49\textwidth}
        \includegraphics[width=\textwidth]{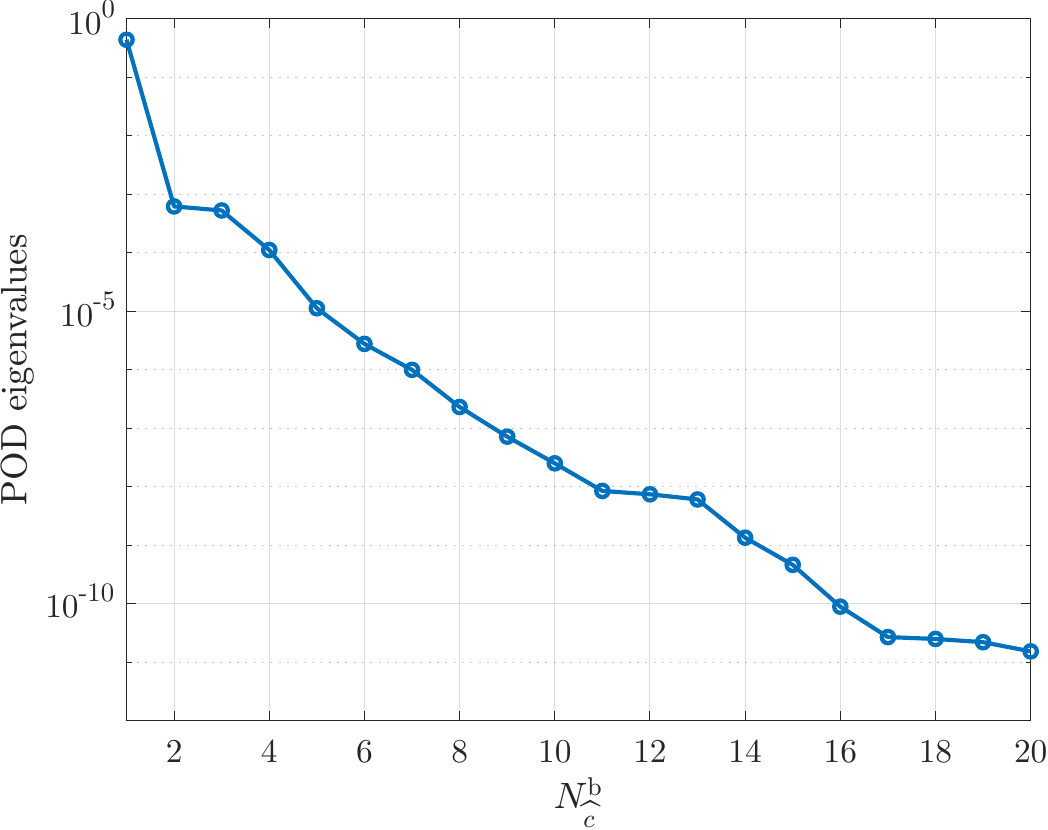}
     \vspace{-1.5em} \caption{Rod component}\label{fig:podrod}
    \end{subfigure}
    \begin{subfigure}{0.49\textwidth}
        \includegraphics[width=\textwidth]{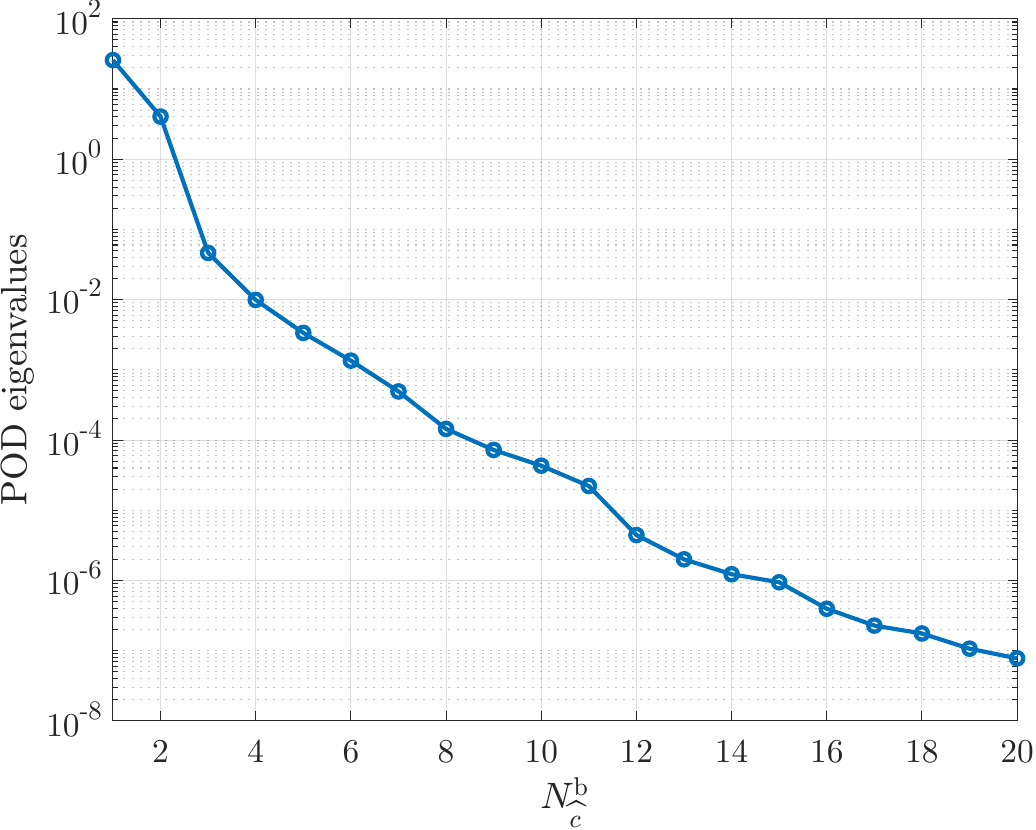}
        \vspace{-1.5em}\caption{Bracket component}\label{fig:podangle}
    \end{subfigure}
    \vspace{1em}
    \begin{subfigure}{0.49\textwidth}
        \includegraphics[width=\textwidth]{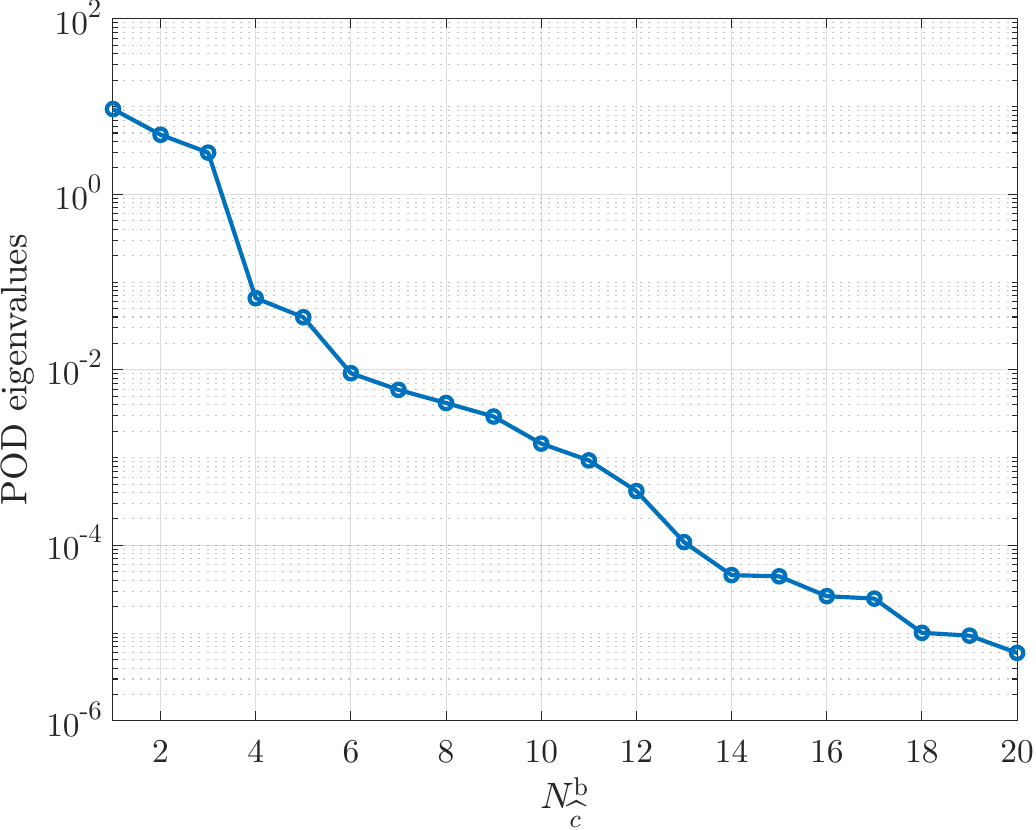}
        \vspace{-1.5em}\caption{Tee component}\label{fig:podtee}
    \end{subfigure}
    \begin{subfigure}{0.49\textwidth}
        \includegraphics[width=\textwidth]{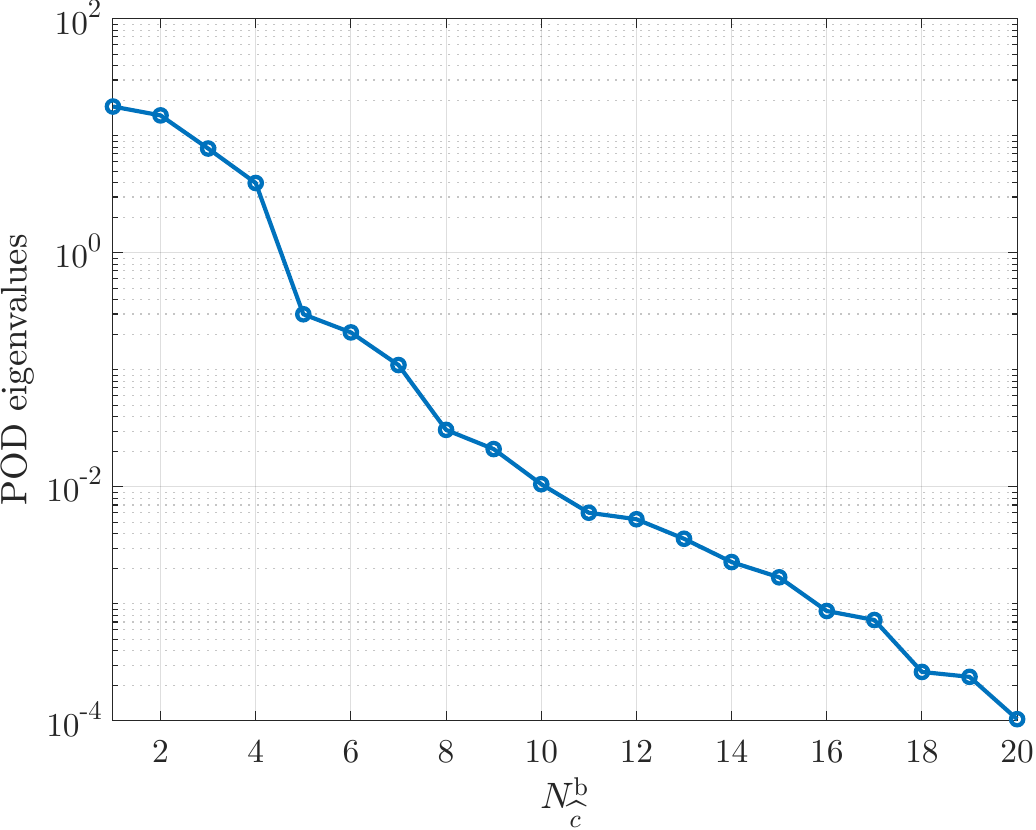}
        \vspace{-1.5em}\caption{Cross component}\label{fig:podcross}
    \end{subfigure}
    \caption{Decay of POD eigenvalues in the RB construction for the bubble space of different archetype components.}
    \label{fig:pods}
\end{figure}

Second, we construct an RB for the bubble space of each archetype component using the POD capturing 99.9\% of the energy (i.e., the sum of POD eigenvalues) of the correlation matrix associated with its bubble snapshot matrix. This results in $N^{\bb}_{\mathrm{rod}} = 3$, $N^{\bb}_{\mathrm{bracket}} = 3$, $N^{\bb}_{\mathrm{tee}} = 6$, and $N^{\bb}_{\mathrm{cross}} = 9$. Figure~\ref{fig:pods} illustrates the decay of POD eigenvalues for each archetype component, showing a rapid decrease in the POD eigenvalues for all components.

\begin{figure}
\centering
    \begin{subfigure}{0.35\textwidth}
        \includegraphics[width=\textwidth]{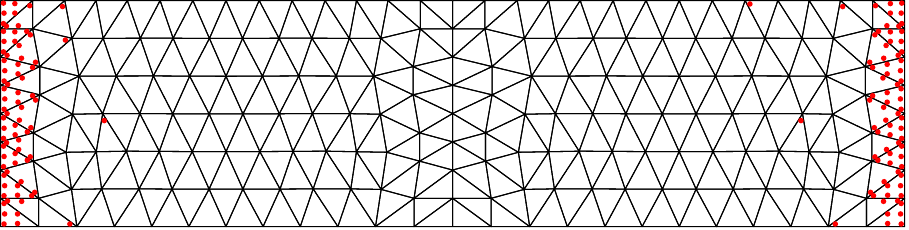}
        \caption{Rod component, $\delta_{\wc} = 10^{2}$}\label{fig:eqprod1}
    \end{subfigure}
    \begin{subfigure}{0.35\textwidth}
        \includegraphics[width=\textwidth]{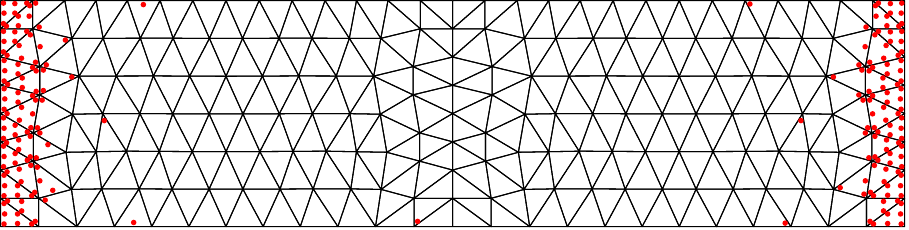}
        \caption{Rod component, $\delta_{\wc} = 1$}\label{fig:eqprod2}
    \end{subfigure}\\[0.5em]
    \begin{subfigure}{0.35\textwidth}
        \includegraphics[width=\textwidth]{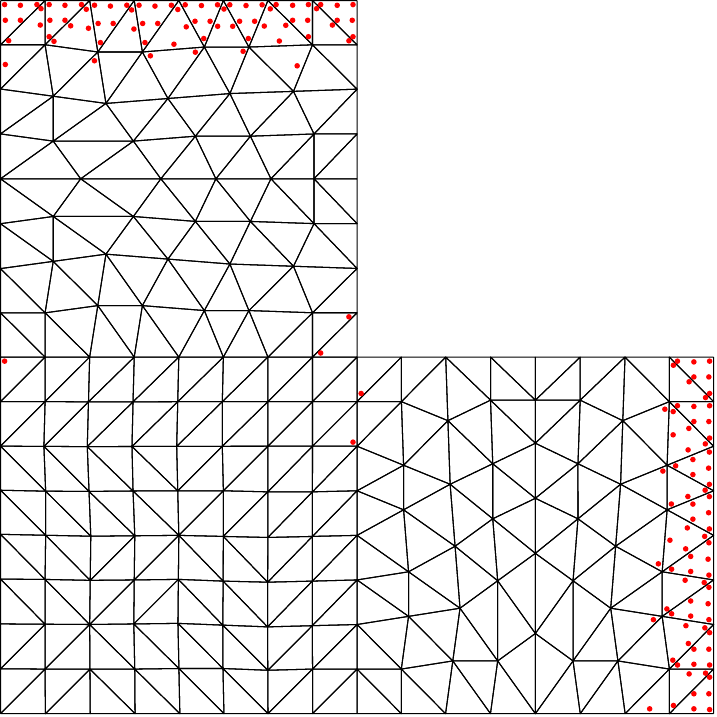}
        \caption{Bracket component, $\delta_{\wc} = 10^{2}$}\label{fig:eqpangle1}
    \end{subfigure}
     \begin{subfigure}{0.35\textwidth}
        \includegraphics[width=\textwidth]{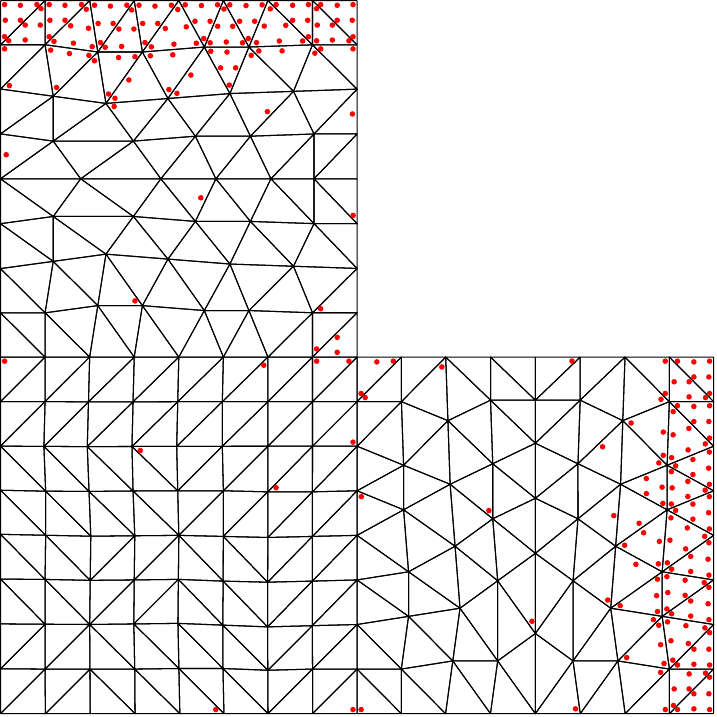}
        \caption{Bracket component, $\delta_{\wc} = 1$}\label{fig:eqpangle2}
    \end{subfigure}\\[0.5em]
    \begin{subfigure}{0.35\textwidth}
        \includegraphics[width=\textwidth]{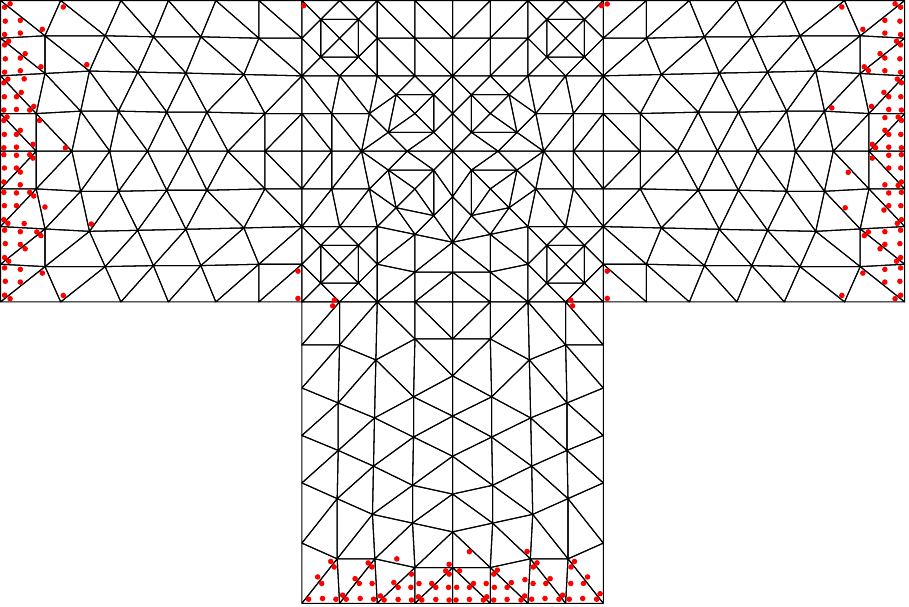}
        \caption{Tee component, $\delta_{\wc} = 10^{2}$}\label{fig:eqptee1}
    \end{subfigure}
    \begin{subfigure}{0.35\textwidth}
        \includegraphics[width=\textwidth]{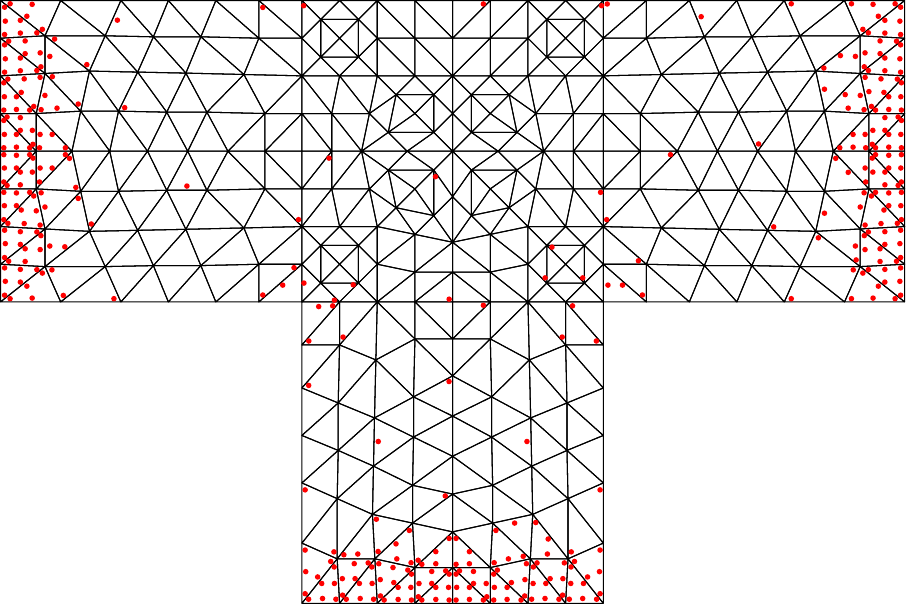}
        \caption{Tee component, $\delta_{\wc} = 1$}\label{fig:eqptee2}
    \end{subfigure}\\[0.5em]
    \begin{subfigure}{0.35\textwidth}
        \includegraphics[width=\textwidth]{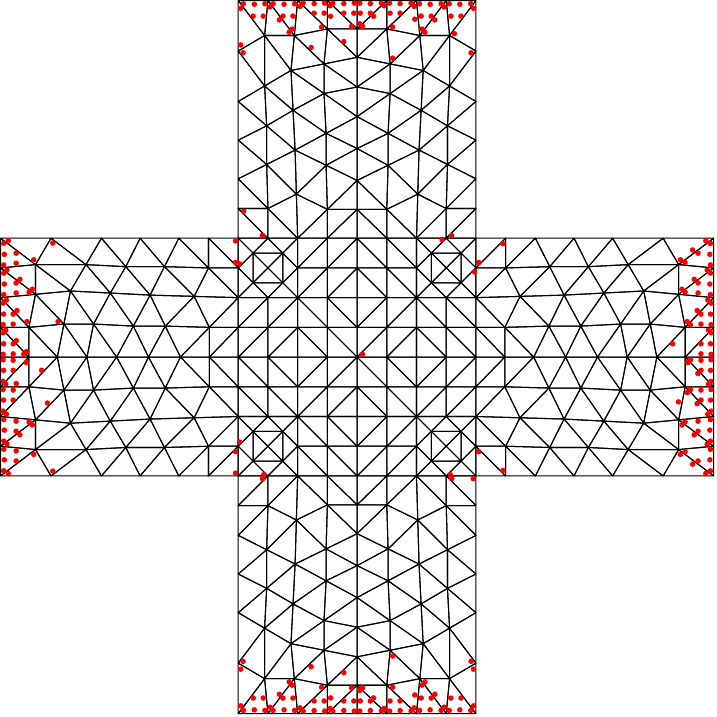}
        \caption{Cross component, $\delta_{\wc} = 10^{2}$}\label{fig:eqpcross1}
    \end{subfigure}
    \begin{subfigure}{0.35\textwidth}
        \includegraphics[width=\textwidth]{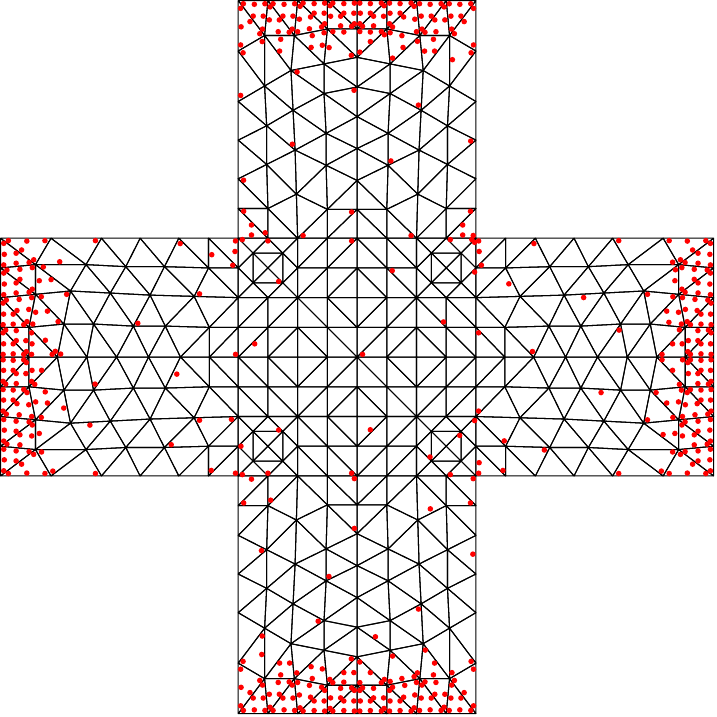}
        \caption{Cross component, $\delta_{\wc} = 1$}\label{fig:eqpcross2}
    \end{subfigure}
    \caption{RQ points of the archetype components for $\delta_\wc = 10^{2}$ and $\delta_\wc = 1$.}
    \label{fig:eqps}
\end{figure}

Third, we follow Algorithm~\ref{alg:emptreqp} to create a bubble RB snapshot set for each archetype component using the data generated in the first step and the RB constructed in the second step. Then, we solve the component-wise hyperreduction problem~\eqref{eq:finaleqp}--\eqref{eq:compconst21} for seven different hyperreduction tolerances $\delta_\wc = \{ 10^{-4}, 10^{-3}, \dots, 10^{2} \}$ to construct a family of RQ rules. Figure~\ref{fig:eqps} shows components' RQ points for hyperreduction tolerances $\delta_\wc = 10^{2}$ and $\delta_\wc = 1$. Table~\ref{table:archsummary} summarizes the outcome of offline training for all archetype components. For all components, the number of bubble degrees of freedom is significantly reduced (i.e., $N^\bb_{\wc} \ll \mathcal{N}^\bb_{\wc}$), and the number of RQ points increases as the hyperreduction tolerance $\delta_{\wc}$ tightens.  

\begin{table}
\caption{\label{table:archsummary}Outcome of offline training for all archetype components.}
\begin{center}
\begin{tabular}{lcccc}
\toprule[0.1em]
{Component} & {Rod} & {Bracket} & {Tee} & {Cross}\\[0.2em]
\hline
$\mathcal{N}^{\bb}_{\wc}$ & 691 & 703 & 1026 & 1165 \\[0.2em]
$Q_{\wc}$ & 1968 & 2016 & 3024 & 3456 \\[0.0em]
\hline \\[-1.1em]
$N^{\bb}_{\wc}$ & 3 & 3 & 6 & 9 \\[0.2em]
$\wtilde{Q}^r_{\wc} \: (\delta_{\wc} = 10^{2})$ & 147 & 156 & 230 & 296 \\[0.2em]
$\wtilde{Q}^r_{\wc} \: (\delta_{\wc} = 10)$ & 183 & 198 & 317 & 420 \\[0.2em]
$\wtilde{Q}^r_{\wc} \: (\delta_{\wc} = 1)$ & 203 & 265 & 391 & 563 \\[0.2em]
$\wtilde{Q}^r_{\wc} \: (\delta_{\wc} = 10^{-1})$ & 287 & 322 & 516 & 739 \\[0.2em]
$\wtilde{Q}^r_{\wc} \: (\delta_{\wc} = 10^{-2})$ & 347 & 375 & 637 & 956 \\[0.2em]
$\wtilde{Q}^r_{\wc} \: (\delta_{\wc} = 10^{-3})$ & 419 & 488 & 846 & 1233\\[0.2em]
$\wtilde{Q}^r_{\wc} \: (\delta_{\wc} = 10^{-4})$ & 482 & 599 & 1031 & 1613\\[0.2em]
\bottomrule
\end{tabular}
\end{center}
\end{table}

\subsection{Thermal fin systems}
\label{subsec:example}
\begin{figure}
  \centering
   \begin{subfigure}{0.49\textwidth}
     \includegraphics[width=1.0\textwidth]{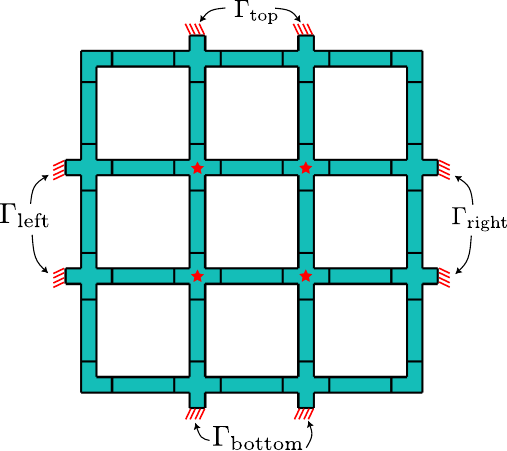}
     \caption{system composition \label{fig:examplefins}}
   \end{subfigure}
   \begin{subfigure}{0.49\textwidth}
     \includegraphics[width=1.0\textwidth]{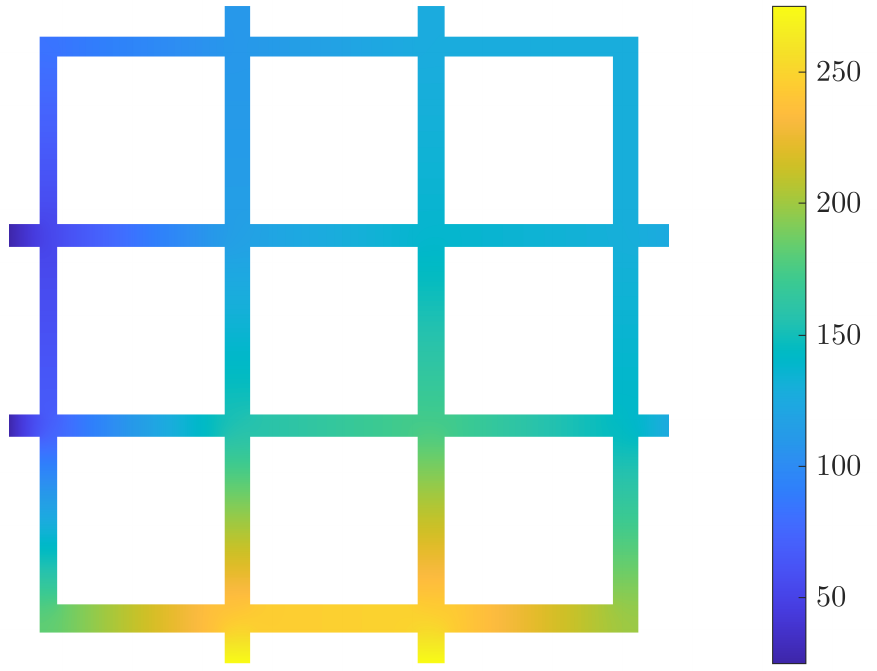}
     \caption{example temperature distribution \label{fig:test1tr}}
   \end{subfigure}
   \caption{A $3 \times 3$ fin system. In (a), red stars mark the components with a volumetric source term.}
\end{figure}

We now examine the performance of the HRBE method on a family of thermal fin systems made of instances of rod, bracket, and cross components from the library. An example of a $3 \times 3$ fin system is shown in Figure~\ref{fig:examplefins}.   We characterize the topology of the fin systems by their number of rod components along horizontal and vertical directions. We consider only the cases where the number of horizontal and vertical rods are identical.

We assume the interior cross components of the fins are subject to a volumetric heat source. Furthermore, we assume the length of the rods along all directions are identical. Additionally, we assume the horizontal and vertical thicknesses vary independently. Hence, an $N_{\fin} \times N_{\fin}$ fin system has $N_{\comp} = (3N_{\fin} + 1) \times (N_{\fin}+1)$ instantiated components, $N_{\fin} + 1$ thickness variables along the horizontal direction, $N_{\fin} + 1$ thickness variables along the vertical direction, 1 length variable associated with rod components, and $(N_{\fin} - 1)^2$ physical variables for volumetric source terms. Therefore, in total, an $N_{\fin} \times N_{\fin}$ fin system is parameterized by $N_{\fin}^2 + 4$ variables, making the problem parametrically high-dimensional even for $N_{\fin} = 2$. Fin systems are subject to four Dirichlet boundary conditions: $u_{\mathrm{left}} = 25$~K on $\Gamma_{\mathrm{left}}$, $u_{\mathrm{right}} = 125$~K on $\Gamma_{\mathrm{right}}$, $u_{\mathrm{bottom}} = 275$~K on $\Gamma_{\mathrm{bottom}}$, and $u_{\mathrm{top}} = 100$~K on $\Gamma_{\mathrm{top}}$. Figure~\ref{fig:test1tr} shows the truth temperature distribution for one instantiation of the $N_{\fin} \times N_{\fin}$ fin system for $N_{\fin} = 3$.

\subsection{Numerical results using prescribed hyperreduction tolerances}
\label{subsubsec:part_1}
\begin{figure}
  \centering
  \begin{subfigure}{0.49\textwidth}
    \includegraphics[width=1.0\textwidth]{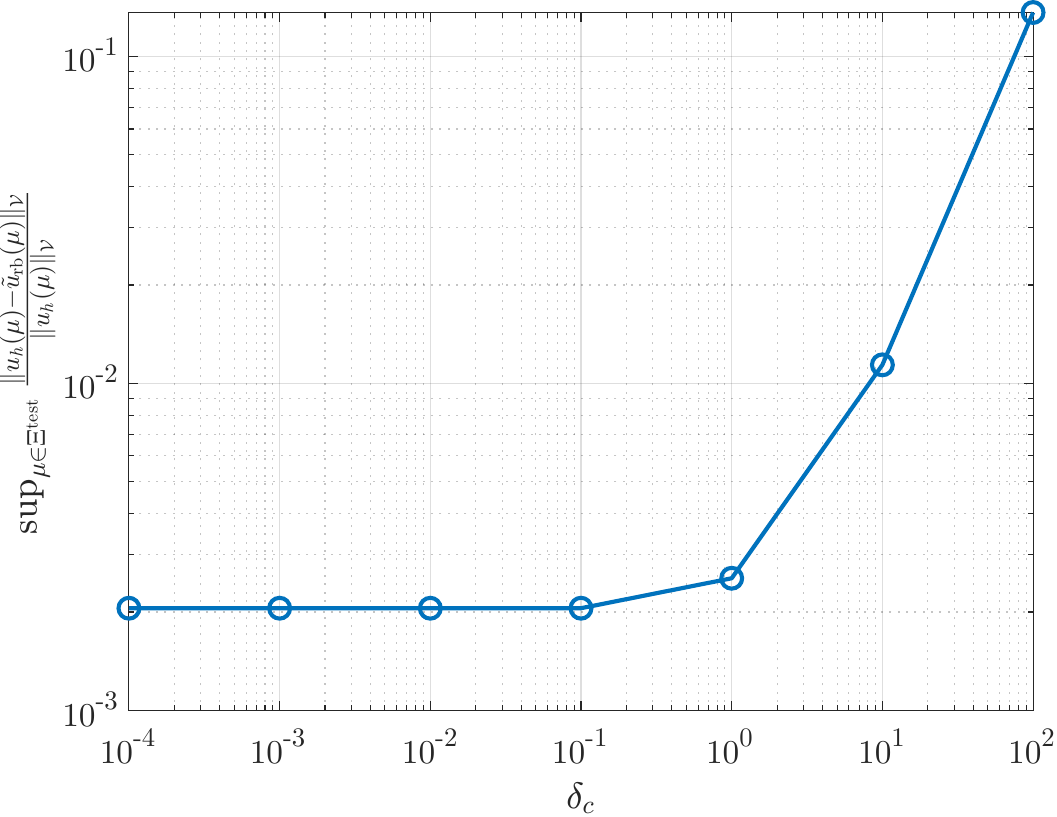}
    \caption{truth vs HRBE \label{fig:hrbetr_erros}}
  \end{subfigure}
  \begin{subfigure}{0.49\textwidth}
    \includegraphics[width=1.0\textwidth]{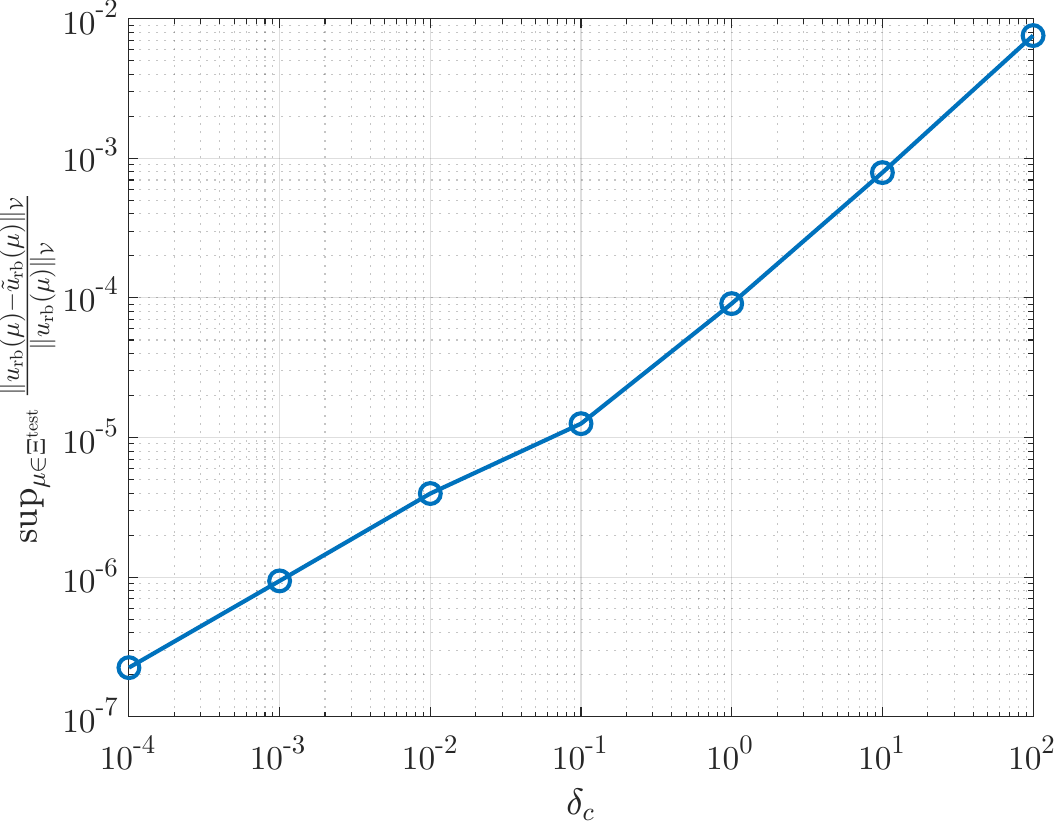}
    \caption{truth-quadrature RB vs HRBE \label{fig:hrberb_erros}}
  \end{subfigure}
  \caption{Maximum relative $H^1(\Omega)$-norm error in the HRBE solution with respect to the truth and RB solutions for different hyperreduction tolerances for the $3 \times 3$ fin over $|\Xi^{\rm test}|=5$ test cases.}
\end{figure}

We first study the behavior of the HRBE method on the $3 \times 3$ fin system using prescribed hyperreduction tolerances $\delta_{c}$ $\forall c \in \sys$; i.e., the same $\delta_c$ is prescribed to all components without using the adaptive algorithm (Algorithm~\ref{alg:sing}). Figure~\ref{fig:hrbetr_erros} shows the relative $H^1(\Omega)$-norm error between the truth and HRBE solutions for different hyperreduction tolerances. To assess the generality of the formulation, we report the maximum error over a test configuration set $\Xi^{\test} \subset \mathcal{D}$, which comprises $|\Xi^{\test}|= 5$ test configurations that results from parameter values randomly selected from a uniform distribution over their corresponding training range. As expected, the error decreases as the hyperreduction tolerances are reduced. A maximum error of $1.363 \times 10^{-1}$ is observed for $\delta_{c} = 10^{2}$, which sharply decreases to $2.536 \times 10^{-3}$ for $\delta_{c} = 1$. The HRBE error eventually plateaus and approaches that of (truth-quadrature) RB without hyperreduction. (We recall that, in this work, the RB is fixed independent of the hyperreduction tolerance for each component, and hence the error between the truth and RB solutions (and in turn the HRBE solutions) is not adaptively controlled; cf.~Remark~\ref{rem:no_adaptive_rb}.)

Figure~\ref{fig:hrberb_erros} shows the maximum relative $H^1(\Omega)$-norm error between (truth-quadrature) RB and HRBE solutions over the $|\Xi^{\test}|=5$ test cases for different $\delta_{c}$ values. As anticipated, the error decreases with hyperreduction tolerances and, hence, when more RQ points are used. More quantitatively, the BRR error bound~\eqref{eq:brrerr} suggests that 
\begin{equation}\nonumber
\left\| {u}_{\rb}(\mu) - \wtilde{{u}}_{\rb}(\mu) \right\|_{\mathcal{V}} \leq \left\| \mathbf{u}_{\rb}(\mu) - \wtilde{\mathbf{u}}_{\rb}(\mu) \right\|_{2} \sqrt{\lambda_\mx}\leq \frac{2\sqrt{\lambda_\mx}}{\sigma_{\min} \left(\wtilde{\mathbf{J}}_{\rb}(\mathbf{u}_{\rb}(\mu); \mu)\right)}\left\| \wtilde{\mathbf{R}}_{\rb}(\mathbf{u}_{\rb}(\mu); \mu) \right\|_2.
\end{equation}
Then, assuming the residual-tolerance condition~\eqref{eq:const1} holds for $\wtilde{\mathbf{u}}_{\rb}(\mu)$ at the system level, we conclude that
\begin{equation}\nonumber
\left\| {u}_{\rb}(\mu) - \wtilde{{u}}_{\rb}(\mu) \right\|_{\mathcal{V}} \leq  \frac{2 \delta_{c} \sqrt{\lambda_\mx} \sum_{c \in \sys} \sqrt{N_{M(c)}}}{\sigma_{\min} \left(\wtilde{\mathbf{J}}_{\rb}(\mathbf{u}_{\rb}(\mu); \mu)\right)}
\end{equation}
since the same $\delta_c$ is applied for all components. 
Hence, if $\sigma_{\min} \left(\wtilde{\mathbf{J}}_{\rb}(\mathbf{u}_{\rb}(\mu); \mu)\right)$ remains approximately constant for different $\delta_{c}$ values, then we expect the error to vary linearly with $\delta_{c}$. This is precisely what we observe in Figure~\ref{fig:hrberb_erros}. The values of $\sigma_{\min} \left(\wtilde{\mathbf{J}}_{\rb}(\mathbf{u}_{\rb}(\mu); \mu)\right)$ reported in Table~\ref{table:sigmas} confirm that the minimum singular value is approximately constant for $\delta_c \leq 10$.

\begin{table}
    \caption{\label{table:sigmas}Value of $\sigma_{\min} \left(\wtilde{\mathbf{J}}_{\rb}(\mathbf{u}_{\rb})\right) \equiv \sigma_{\min} \left(\wtilde{\mathbf{J}}_{\rb}(\mathbf{u}_{\rb}(\mu); \mu)\right)$ for different hyperreduction tolerances for the $3 \times 3$ fin.}
    \begin{center}
    \begin{tabular}{lccccccc}
    \toprule[0.1em]
    \begin{tabular}{@{}l} Hyperreduction\\{tolerances} \end{tabular} & $\delta_{c} = 10^{2}$ & $\delta_{c} = 10$ & $\delta_{c} = 1$ & $\delta_{c} = 10^{-1}$ & $\delta_{c} = 10^{-2}$ & $\delta_{c} = 10^{-3}$ & $\delta_{c} = 10^{-4}$ \\[0.2em]
    $\sigma_{\min} \left(\wtilde{\mathbf{J}}_{\rb}(\mathbf{u}_{\rb})\right)$ &  $2.697$ & $2.927$ & $2.932$ & $2.933$ & $2.933$ & $2.933$ & $2.933$ \\
    \bottomrule
    \end{tabular}
    \end{center}
\end{table}

We now study the behavior of the minimum singular value $\sigma_{\mn}(\wtilde{\mathbf{J}}_{\rb}\left(\wtilde{\mathbf{u}}_{\rb}(\mu);\mu) \right)$, which plays an important role in the adaptive RQ selection in Algorithm~\ref{alg:sing}. To develop the algorithm, we posited, based on Lemma~\ref{lemma:singularpert}, that $\sigma_{\mn}(\wtilde{\mathbf{J}}_{\rb}\left(\wtilde{\mathbf{u}}_{\rb}(\mu);\mu) \right)$ would provide a reliable approximation for $\sigma_{\mn}({\mathbf{J}}_{\rb}\left(\wtilde{\mathbf{u}}_{\rb}(\mu);\mu) \right)$. Table~\ref{table:hrbe_sigmas} shows the maximum relative error between these two values $\forall \mu \in \Xi^{\test}$ for different hyperreduction tolerances~$\delta_c$. We note that even for the highest $\delta_{c}$ the error between the singular values is less than $10\%$ and the difference quickly decreases as $\delta_{c}$ is reduced. Consequently, in practice, as Algorithm~\ref{alg:sing} iterates toward smaller hyperreduction tolerances, this error becomes increasingly insignificant. Additionally, the applied factor of $0.5$ due to~\eqref{eq:sigmaest} further mitigates the possibility that $0.5 \sigma_{\mn}(\wtilde{\mathbf{J}}_{\rb}\left(\wtilde{\mathbf{u}}_{\rb}(\mu);\mu) \right)$ does not provide a lower bound of $\sigma_{\mn}({\mathbf{J}}_{\rb}\left({\mathbf{u}}_{\rb}(\mu);\mu) \right)$ in Algorithm~\ref{alg:sing}.

\begin{table}
\caption{\label{table:hrbe_sigmas}Relative error between $\sigma_{\mn}(\wtilde{\mathbf{J}}_{\rb}\left(\wtilde{\mathbf{u}}_{\rb}(\mu);\mu) \right)$ and $\sigma_{\mn}({\mathbf{J}}_{\rb}\left(\wtilde{\mathbf{u}}_{\rb}(\mu);\mu) \right)$ for different hyperreduction tolerances for the $3 \times 3$ fin over $|\Xi^{\rm test}|=5$ test cases.}
\begin{center}
\begin{tabular}{lc}
\toprule[0.1em]
\begin{tabular}{@{}l} Hyperreduction\\{tolerances} \end{tabular} & $\displaystyle \sup_{\mu \in \Xi^{\test}} \frac{|\sigma_{\mn}({\mathbf{J}}_{\rb}\left(\wtilde{\mathbf{u}}_{\rb}(\mu);\mu) \right) - \sigma_{\mn}(\wtilde{\mathbf{J}}_{\rb}\left(\wtilde{\mathbf{u}}_{\rb}(\mu);\mu) \right)|}{\sigma_{\mn}({\mathbf{J}}_{\rb}\left(\wtilde{\mathbf{u}}_{\rb}(\mu);\mu) \right)}$ \\[0.2em]
\hline
$\delta_{c} = 10^{2}$ &  $9.009 \times 10^{-2}$ \\[0.2em]
$\delta_{c} = 10$ & $3.017 \times 10^{-3}$ \\[0.2em]
$\delta_{c} = 1$ & $1.507 \times 10^{-4}$ \\[0.2em]
$\delta_{c} = 10^{-1}$ & $1.057 \times 10^{-5}$ \\[0.2em]
$\delta_{c} = 10^{-2}$ &  $1.023 \times 10^{-6}$ \\[0.2em]
$\delta_{c} = 10^{-3}$ & $6.818 \times 10^{-7}$ \\[0.2em]
$\delta_{c} = 10^{-4}$ & $4.973 \times 10^{-7}$ \\
\bottomrule
\end{tabular}
\end{center}
\end{table}

\begin{figure}
\centering
\includegraphics[width=0.49\textwidth]{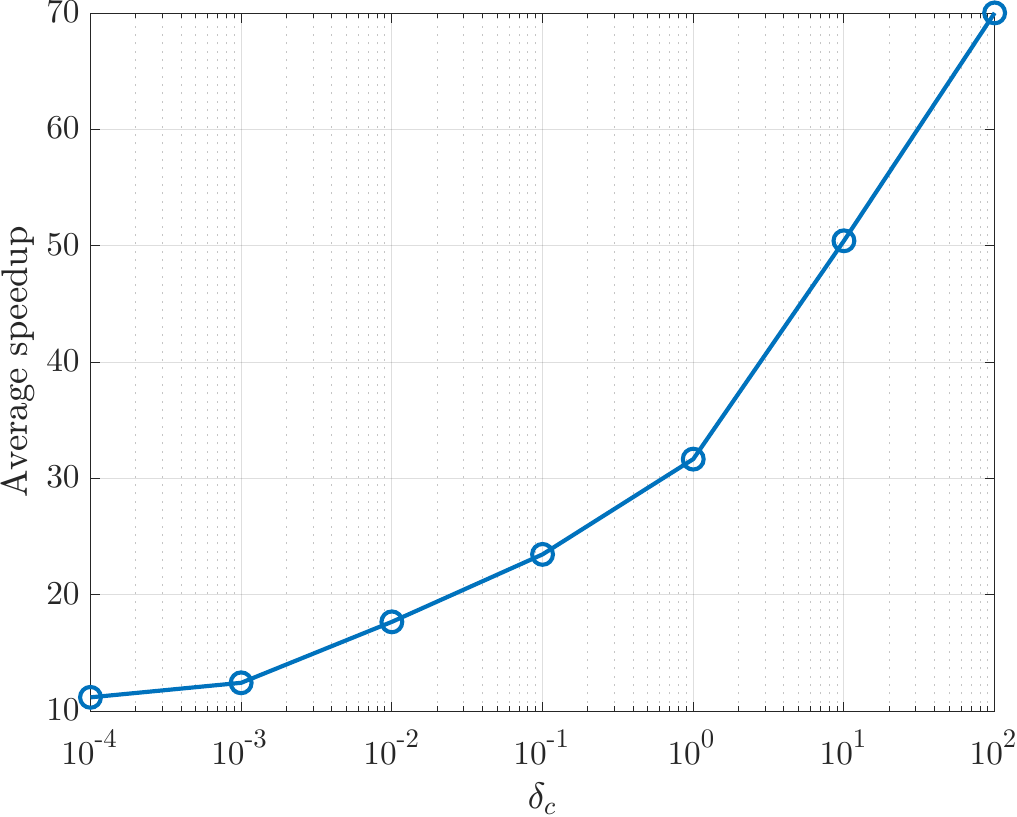}
        \caption{Average speedup in wall-clock time relative to solving the truth problem for different hyperreduction tolerances for the $3 \times 3$ fin across $|\Xi^{\rm test}|=5$ test cases.}\label{fig:hrbe_speeds}
\end{figure}

Figure~\ref{fig:hrbe_speeds} shows the average speedup in wall-clock time relative to solving the truth problem across the five test configurations for various hyperreduction tolerances. Specifically, an average speedup of around 70 times is observed for $\delta_{c} = 10^{2}$, reducing to about 11 times for $\delta_{c} = 10^{-4}$. While the difference in speedups might encourage the use of looser hyperreduction tolerances, it is crucial to consider the trade-off in accuracy. We recall that Figures~\ref{fig:hrbetr_erros} and \ref{fig:hrberb_erros} show that the errors for $\delta_{c} = 10^{2}$ is significantly higher than the errors for the RQ rules associated with tighter tolerances. Conversely, opting for the strictest tolerance yields the most accurate HRBE solution, but the speedup is not as substantial compared to using looser tolerances. (We recall that, in this work, we do not consider port reduction, and hence the speedup achieved by the HRBE method is $\mathcal{O}(10)$--$\mathcal{O}(100)$ and not $\calO(1000)$ as achieved by port-reduced RBEs for linear problems~\cite{eftang2014port}; cf.~Remark~\ref{rem:no_port_red}.)

\subsection{Numerical results using the adaptive RQ selection algorithm}
\label{subsubsec:part_2}
To effectively navigate the trade-off between speedup and accuracy and select the RQ rules satisfying a desired error between RB and HRBE solutions, we now apply Algorithm~\ref{alg:sing} for the relative error target $\epsilon = 0.01$. The algorithm finds the RQ rules corresponding to different hyperreduction tolerances for different components in the $3 \times 3$ thermal fin system, although the same tolerance is used for the components instantiated from the same archetype component. Convergence is reached in merely two iterations in all parameter test configurations. The maximum relative $H^1(\Omega)$-norm errors in the HRBE solutions relative to the truth and (truth-quadrature) RB solutions are $7.521 \times 10^{-3}$ and $7.226 \times 10^{-3}$, respectively. As desired, the adaptive RQ selection algorithm meets the target \emph{system-level} hyperreduction error tolerance of~$10^{-2}$. The HRBE method provides an average computational speedup of 42 relative to solving the truth problem.

To further assess the performance of the adaptive RQ selection algorithm across a range of fin system sizes, we apply Algorithm~\ref{alg:sing}, with the relative error target $\epsilon = 0.01$, to $N_{\mathrm{fin}} \times N_{\mathrm{fin}}$ fin systems for $N_{\mathrm{fin}} \in  \{2, \cdots,8\}$. For each fin system, we form $|\Xi^{\test}| = 5$ test configurations similar to those for the $3 \times 3$ fin system described earlier. For all fin systems, the algorithm achieves convergence within two iterations. Table~\ref{table:eqplevels} shows the maximum relative $H^1(\Omega)$-norm errors across the test configurations between (i) truth and RB solutions, (ii) truth and HRBE solutions, and (iii) RB and HRBE solutions. The target error between the RB and HRBE solutions is achieved for all fin systems. The effectivity, defined as $\epsilon$ divided by the actual maximum relative error, ranges from a minimum of $1.315$ for $N_{\fin} = 2$ to a maximum of $33.602$ for $N_{\fin} = 8$. The sharpness of the error bound between the RB and HRBE solutions deteriorates as $N_{\mathrm{fin}}$ increases. We suspect that this is due to bounding $\lr{ \cdot }_{\infty}$ and $\lr{ \cdot }_{\mx}$ in Proposition~\ref{prop:suffcond} by $\sqrt{N_{M(c)}} \lr{ \cdot }_{2}$ and $N_{M(c)} \lr{\cdot }_{2}$ $\forall c \in \sys$, respectively. As the number of components in the system increases, $\bar{\alpha}$ in Proposition~\ref{prop:suffcond} provides a more pessimistic upper bound for the error at the system level. Table~\ref{table:eqplevels} also shows that for all fin systems the error between the truth and RB solutions is relatively close to the error between the truth and HRBE solutions, underscoring the effectiveness of the adaptive RQ selection (Algorithm~\ref{alg:sing}) as well as the component-wise hyperreduction training routine. 

\begin{table}
  \caption{\label{table:eqplevels}Relative $H^1(\Omega)$-norm error between (i) truth and RB solutions, (ii) truth and HRBE solutions, and (iii) RB and HRBE solutions for $N_{\mathrm{fin}} \times N_{\mathrm{fin}}$ fins using $\epsilon = 0.01$ in Algorithm~\ref{alg:sing} over their five test cases.}
\begin{center}
\begin{tabular}{lccc}
\toprule[0.1em]
$N_{\mathrm{fin}}$ & $\sup_{\mu \in \Xi^{\test}} \frac{\|u_h(\mu) - {u}_{\rb}(\mu) \|_{\mathcal{V}}}{\|u_h(\mu) \|_{\mathcal{V}}}$ & $\sup_{\mu \in \Xi^{\test}} \frac{\|u_h(\mu) - \wtilde{u}_{\rb}(\mu) \|_{\mathcal{V}}}{\|u_h(\mu) \|_{\mathcal{V}}}$ & $\sup_{\mu \in \Xi^{\test}} \frac{\|u_\rb(\mu) - \wtilde{u}_{\rb}(\mu) \|_{\mathcal{V}}}{\|u_\rb(\mu) \|_{\mathcal{V}}}$ \\[0.2em]
\hline
$2$ & $4.183 \times 10^{-3}$ & $7.891 \times 10^{-3}$ & $7.605 \times 10^{-3}$ \\[0.2em]
$3$ & $2.054 \times 10^{-3}$ & $7.521 \times 10^{-3}$ & $7.226 \times 10^{-3}$ \\[0.2em]
$4$ & $1.557 \times 10^{-3}$ & $6.134 \times 10^{-3}$ & $5.921 \times 10^{-3}$ \\[0.2em]
$5$ & $1.626 \times 10^{-3}$ & $4.958 \times 10^{-3}$ & $4.832 \times 10^{-3}$ \\[0.2em]
$6$ & $1.283 \times 10^{-3}$ & $1.341 \times 10^{-3}$ & $3.801 \times 10^{-4}$ \\[0.2em]
$7$ & $1.001 \times 10^{-3}$ & $1.057 \times 10^{-3}$ & $3.377 \times 10^{-4}$ \\[0.2em]
$8$ & $8.882 \times 10^{-4}$ & $9.389 \times 10^{-4}$ & $2.976 \times 10^{-4}$ \\
\bottomrule
\end{tabular}
\end{center}
\end{table}

Finally, Table~\ref{table:sigmafins} presents the relative error between $\sigma_{\mn}(\wtilde{\mathbf{J}}_{\rb}\left(\wtilde{\mathbf{u}}_{\rb}(\mu);\mu) \right)$ and $\sigma_{\mn}({\mathbf{J}}_{\rb}\left(\wtilde{\mathbf{u}}_{\rb}(\mu);\mu) \right)$ across various fin system sizes, which again supports the validity of using $\sigma_{\mn}(\wtilde{\mathbf{J}}_{\rb}\left(\wtilde{\mathbf{u}}_{\rb}(\mu);\mu) \right)$ to approximate $\sigma_{\mn}({\mathbf{J}}_{\rb}\left(\wtilde{\mathbf{u}}_{\rb}(\mu);\mu) \right)$ in Algorithm~\ref{alg:sing}.

\begin{table}[!tb]
    \caption{\label{table:sigmafins}Relative error between $\sigma_{\mn}(\wtilde{\mathbf{J}}_{\rb}\left(\wtilde{\mathbf{u}}_{\rb}(\mu);\mu) \right)$ and $\sigma_{\mn}({\mathbf{J}}_{\rb}\left(\wtilde{\mathbf{u}}_{\rb}(\mu);\mu) \right)$ for $N_{\mathrm{fin}} \times N_{\mathrm{fin}}$ fins over their five test cases.}
    \begin{center}
    \begin{tabular}{lc}
    \toprule[0.1em]
    $N_{\mathrm{fin}}$ & $\displaystyle \sup_{\mu \in \Xi^{\test}} \frac{|\sigma_{\mn}({\mathbf{J}}_{\rb}\left(\wtilde{\mathbf{u}}_{\rb}(\mu);\mu) \right) - \sigma_{\mn}(\wtilde{\mathbf{J}}_{\rb}\left(\wtilde{\mathbf{u}}_{\rb}(\mu);\mu) \right)|}{\sigma_{\mn}({\mathbf{J}}_{\rb}\left(\wtilde{\mathbf{u}}_{\rb}(\mu);\mu) \right)}$ \\[0.2em]
    \hline
    $2$ & $3.478 \times 10^{-3}$ \\[0.2em]
    $3$ & $2.395 \times 10^{-3}$ \\[0.2em]
    $4$ & $2.159 \times 10^{-3}$ \\[0.2em]
    $5$ & $1.484 \times 10^{-3}$ \\[0.2em]
    $6$ & $4.221 \times 10^{-5}$ \\[0.2em]
    $7$ & $4.955 \times 10^{-5}$ \\[0.2em]
    $8$ & $6.616 \times 10^{-5}$ \\
    \bottomrule
    \end{tabular}
    \end{center}
\end{table}

\section{Conclusion}
\label{sec:conclusion}
In this work, we have developed an HRBE method for reduced-order modeling of component-based systems governed by general parameterized nonlinear PDEs. The proposed method is capable of accommodating global nonlinearities across the entire domain. The method constructs a library of archetype components during the offline phase through component-wise RB construction and hyperreduction. Then, in the online phase, these pretrained components are reused to rapidly create a reduced model for any system configuration instantiated from the archetype components in the library. This divide-and-conquer strategy circumvents the need for repeated offline training for new system configurations and enables the reduced-order modeling of problems with numerous parameters. Additionally, it facilitates the model reduction of large-scale problems by sidestepping the generation of global solution snapshots associated with large assembled systems in the offline phase.

The proposed HRBE method is characterized by several key features. First, we have formulated a component-wise extension of the EQP~\cite{Patera_2017_EQP,yano2019lp} to systematically construct a library of hyperreduced components, each of which meets the specified hyperreduction tolerance. Second, we have appealed to the BRR theorem to develop an actionable error estimate for component-based systems, which relates component-wise hyperreduction residual to the system-level error. Third, we have developed an online-efficient estimate of the minimum singular value of the system-level Jacobian, which plays a crucial role in the BRR theory. Finally, building on the aforementioned multi-fidelity archetype component library, the actionable error estimate, and the minimum singular value estimate, we have developed an adaptive RQ selection procedure, such that the hyperreduction error in the online-assembled system meets the user-prescribed system-level error tolerance. 

We evaluated the effectiveness of our HRBE method through its application to two-dimensional nonlinear thermal fin systems, which are composed from a library consisting of four distinct types of archetype components. Across different fin systems, we demonstrated that the HRBE method consistently delivers accurate results and computational reduction, achieving roughly 45$\times$ speedups with errors around 1\% or less. Moreover, the online-efficient minimum singular value estimate for the system's RB Jacobian proved accurate in the fin systems studied.

There exist several potential opportunities to extend the current work. First is the development of a port-reduced version of the HRBE method (cf.~Remark~\ref{rem:no_port_red}). In systems with many and/or large global ports, the final HRBE problem can still be quite large without port reduction. Hence, model reduction of the ports could lead to additional computational savings in the online phase, a concept already explored for linear problems (e.g.,~\cite{eftang2013port, eftang2014port, smetana2015new}). Second is the development of an online-efficient system-level a posteriori error estimates, which is another area that has been explored for linear problems. Third, building on the a posteriori error estimate, Algorithm~\ref{alg:sing} may be extended to effect adaptive selection of both RB and RQ in the online phase (cf.~Remark~\ref{rem:no_adaptive_rb}). Lastly, the current work could be expanded to accommodate time-dependent nonlinear PDEs. We aim to explore these potential extensions in our future research.

\section*{Acknowledgment}
We would like to thank Prof.~Anthony Patera (MIT) for the initial discussions of component-based model reduction for nonlinear systems.

\newpage
\begin{appendices}
\newpage
\section{Explicit expressions of the algebraic RB and hyperreduced RB residuals and Jacobians of instantiated components}
\label{app:resandjacrb}

We assume in each archetype component $\wc \in \lib$, the generalized coordinates of any $w_{\rb,c} \in \mathcal{V}_{\rb,c}$ are arranged in such a way that the DoF associated with $\widehat{\mathcal{V}}_{\rb,\wc}^\bb$ are assigned to the first $N_{\wc}^\bb$ indices, followed by $\mathcal{N}_{\wc}^1$ indices corresponding to the first local port's DoF and so forth. As such, for the instantiated component $c \in \sys$, we introduce $\{ \Psi_{c,i} \}_{i=1}^{N_{M(c)}}$ for any $\mu_c \in \mathcal{D}_c$ as
\begin{equation}\nonumber
    \Psi_{c,i} = \begin{cases}
    \xi_{c,i}^{\bb} \equiv \widehat{\xi}_{M(c),i}^\bb \circ \mathcal{G}_c^{-1}(\cdot; \mu_c), & \quad i \in \{ 1,\dots,N^\bb_{M(c)}\}, \\[0.4em]
    \psi_{c,i}^{1}, & \quad i \in N^\bb_{M(c)} + \{ 1, \dots, \mathcal{N}_{M(c)}^1 \}, \\[0.4em]
    \: \: \vdots\\[0.4em]
    \psi_{c,i}^{n^\gamma_{M(c)}}, & \quad i \in N^\bb_{M(c)} + \sum_{p = 1}^{n^\gamma_{M(c)} - 1} + \{ 1, \dots, \mathcal{N}_{M(c)}^{n^\gamma_{M(c)}} \}. \\[0.4em]
\end{cases}
\end{equation}
The algebraic RB residual and Jacobian for $c \in \sys$ are then given by
\begin{equation}\label{eq:resbrb}
\begin{alignedat}{3}
\left( {\mathbf{R}}_{\rb,c}({{w}}_{\rb,c}; \mu_c) \right)_i &= \sum_{q=1}^{{Q}_{M(c)}} \what{\rho}_{M(c),q} \:  \widehat{r}_{M(c)} \biggl( \Bigl[ {w}^{\bb}_{\rb,c}+ \sum_{p \in \pset_{M(c)}} w^p_{h,p} \Bigr] \circ {\mathcal{G}}_{c}(\cdot;\mu_c), \Psi_{c,i}  \circ {\mathcal{G}}_{c}(\cdot;\mu_c); \what{x}_{M(c),q},\mu_c \biggr), \\
\left({\mathbf{J}}_{\rb,c}({{w}}_{\rb,c}; \mu_c) \right)_{i,j} &= \sum_{q=1}^{{Q}_{M(c)}} \what{\rho}_{M(c),q} \widehat{r}\:'_{M(c)} \Bigl(\Bigl[ {w}^{\bb}_{\rb,c}+ \sum_{p \in \pset_{M(c)}} w^p_{h,c} \Bigr] \circ {\mathcal{G}}_{c}(\cdot;\mu_c), \Psi_{c,j} \circ {\mathcal{G}}_{c}(\cdot;\mu_c),\Psi_{c,i} \circ {\mathcal{G}}_{c}(\cdot;\mu_c);\\
&\hspace{33em} \what{x}_{M(c),q},\mu_c \Bigr),
\end{alignedat}
\end{equation}
for $1 \leq i,j \leq N_{M(c)}$, any ${w}_{\rb,c} \in \mathcal{V}_{\rb,c}$ and any $\mu_c \in \mathcal{D}_c$. 

Similarly, the algebraic hyperreduced RB residual and Jacobian for $c \in \sys$ are given by
\begin{equation}\label{eq:resbhrbe}
    \begin{alignedat}{3}
    \left( \wtilde{\mathbf{R}}_{\rb,c}({{w}}_{\rb,c}; \mu_c) \right)_i &= \sum_{q=1}^{\wtilde{Q}^r_{M(c)}} \th{\rho}^r_{M(c),q} \:  \widehat{r}_{M(c)} \biggl( \Bigl[ {w}^{\bb}_{\rb,c}+ \sum_{p \in \pset_{M(c)}} w^p_{h,c} \Bigr] \circ {\mathcal{G}}_{c}(\cdot;\mu_c), \Psi_{c,i}  \circ {\mathcal{G}}_{c}(\cdot;\mu_c); \th{x}^r_{M(c),q},\mu_c \biggr), \\
    \left(\wtilde{\mathbf{J}}_{\rb,c}({{w}}_{\rb,c}; \mu_c) \right)_{i,j} &= \sum_{q=1}^{\wtilde{Q}_{M(c)}^r} \th{\rho}^r_{M(c),q} \widehat{r}\:'_{M(c)} \Bigl(\Bigl[ {w}^{\bb}_{\rb,c}+ \sum_{p \in \pset_{M(c)}} w^p_{h,c} \Bigr] \circ {\mathcal{G}}_{c}(\cdot;\mu_c), \Psi_{c,j} \circ {\mathcal{G}}_{c}(\cdot;\mu_c),\Psi_{c,i} \circ {\mathcal{G}}_{c}(\cdot;\mu_c);\\
    &\hspace{33em} \th{x}^r_{M(c),q},\mu_c \Bigr),
\end{alignedat}
\end{equation}
for $1 \leq i,j \leq N_{M(c)}$, any ${w}_{\rb,c} \in \mathcal{V}_{\rb,c}$ and any $\mu_c \in \mathcal{D}_c$. 

\end{appendices}

\newpage
\bibliography{main}
\bibliographystyle{ieeetr}
\end{document}